\documentclass[10pt]{elsarticle}
\usepackage[cp1251]{inputenc}
\usepackage[english]{babel}
\usepackage{amsmath}
\usepackage{amssymb}
\usepackage{amsfonts}
\usepackage{rotating}
\usepackage{graphicx}
\usepackage{floatflt,epsfig}
\usepackage{lineno,hyperref}
\usepackage{enumerate}
\usepackage{colortbl}

\usepackage{lipsum}
\makeatletter
\def\ps@pprintTitle{%
 \let\@oddhead\@empty
 \let\@evenhead\@empty
 \def\@oddfoot{}%
 \let\@evenfoot\@oddfoot}
\makeatother

\usepackage[table]{xcolor}
\usepackage{array,tabularx,tabulary,booktabs}
\usepackage{longtable}
\usepackage{multirow}
\usepackage{wrapfig}
\usepackage{subcaption}
\usepackage{float}
\restylefloat{table}

\newcolumntype{$}{>{\global\let\currentrowstyle\relax}}
\newcolumntype{^}{>{\currentrowstyle}}

\journal{European Journal of Combinatorics}
\newtheorem{lemma}{Lemma}

\newtheorem{theorem}{Theorem}
\newtheorem{corollary}{Corollary}

\newcommand{\proof}{\medskip\noindent{\bf Proof.~}}
\begin{document}
\renewcommand{\abstractname}{Abstract}
\renewcommand{\refname}{References}
\renewcommand{\tablename}{Figure.}
\renewcommand{\arraystretch}{0.9}
\thispagestyle{empty}
\sloppy

\begin{frontmatter}

\title{The smallest strictly Neumaier graph and its generalisations}%\tnoteref{grant}}
\tnotetext[grant]{The second author is supported by the National Science Foundation of China, STCSM (17690740800) and RFBR (17-51-560008).
The third author is supported by RFBR (17-51-560008).}

\author{Rhys~J.~Evans*}%\corref{cor1}}
\cortext[cor1]{Corresponding author}
\address{School of Mathematical Sciences\\Queen Mary University of London\\Mile End Road\\London E14NS, UK\\\nolinkurl{r.evans@qmul.ac.uk}}

\author{Sergey~Goryainov}
\address{School of Mathematical Sciences\\Shanghai Jiao Tong University\\800 Dongchuan Road\\Minhang District, Shanghai, China}
\address{Krasovskii Institute of Mathematics and Mechanics\\S. Kovalevskaja st. 16\\Yekaterinburg, 620990, Russia}
\address{Faculty of Mathematics\\Chelyabinsk State University\\Brat'ev Kashirinyh st. 129\\Chelyabinsk,  454021, Russia\\\nolinkurl{44g@mail.ru}}

\author{Dmitry~Panasenko}
\address{Faculty of Mathematics\\Chelyabinsk State University\\Brat'ev Kashirinyh st. 129\\Chelyabinsk,  454021, Russia\\\nolinkurl{makare95@mail.ru}}

\begin{abstract}
A \emph{regular clique} in a regular graph is a clique such that every vertex outside of the clique is adjacent to the same positive number of vertices inside the clique. We continue the study of regular cliques in edge-regular graphs initiated by A. Neumaier in the 1980s and attracting current interest. We thus define a \emph{Neumaier graph} to be an non-complete edge-regular graph containing a regular clique, and a \emph{strictly Neumaier graph} to be a non-strongly regular Neumaier graph. We first prove some general results on Neumaier graphs and their feasible parameter tuples. We then apply these results to determine the smallest strictly Neumaier graph, which has $16$ vertices. Next we find the parameter tuples for all strictly Neumaier graphs having at most $24$ vertices. Finally, we give two sequences of graphs, each with $i^{\text{th}}$ element a strictly Neumaier graph containing a $2^{i}$-regular clique (where $i$ is a positive integer) and having parameters of an affine polar graph as an edge-regular graph. This answers questions recently posed by G. Greaves and J. Koolen.
\end{abstract}

\begin{keyword}
edge-regular graph; strongly regular graph; regular clique; Neumaier graph; affine polar graph; switching.
\vspace{\baselineskip}
\end{keyword}
\end{frontmatter}
%%%%%%%%%%%%%%%%%%%%%%%%%%%%%%%%%%%%%%%%%%%%%%%%%%%%%%%%%%%%%%%%%%%%%%%%%%%%%%%%%%%%%%%%%%%%%%%%%%%%%%%%%%%%%%%%%%%%%%%%%%%%%%%%%%

\section{Introduction}

A \emph{regular clique}, or more specifically an $m$\emph{-regular clique}, in a regular graph $\Gamma$ is a clique $S$ such that every vertex of $\Gamma$ not in $S$ is adjacent to the same positive number $m$ of vertices of $S$. A regular clique can be equivalently viewed as a clique which is a part of an equitable 2-partition (see \cite{BH12,GR01}), or a completely regular code of radius 1 (see \cite{N92} and \cite[p. 345]{BCN89}). It is well known that a clique in a strongly regular graph is regular if and only if it is a Delsarte clique
(see \cite{BHK07}; \cite[Proposition 1.3.2(ii)]{BCN89}; \cite[Proposition 4.4.6]{BCN89}).
\par
In the early 1980s, A. Neumaier \cite{N81} studied regular cliques in edge-regular graphs, and a certain class of designs whose point graphs are strongly regular and contain regular cliques. He then posed the problem of whether there exists a non-complete, edge-regular, non-strongly regular graph containing a regular clique. We thus define a \emph{Neumaier graph} to be a non-complete edge-regular graph containing a regular clique and define a \emph{strictly Neumaier graph} to be a non-strongly regular Neumaier graph. (This is analogous to the definitions of Deza graphs and strictly Deza graphs \cite{EFHHH99}.)
\par
Informed about the problem by L. Soicher in 2015, G. Greaves and J. Koolen then gave an answer by constructing an infinite family of strictly Neumaier graphs \cite{GK18}. A. Gavrilyuk and S. Goryainov then searched for examples in a collection of known Cayley-Deza graphs \cite{GS14}, leading to the discovery of four more strictly Neumaier graphs.
\par
Further to a discussion with Koolen, Goryainov and his student D. Panasenko found the smallest strictly Neumaier graph, using methods similar to some of their work on Deza graphs. At roughly the same time, Evans found the smallest strictly Neumaier graph in a collection of vertex-transitive edge-regular graphs which had been provided by G. Royle \cite{HR17}. Subsequent communications led to the collaboration found in the current paper.
\par
In their paper \cite{GK18}, Greaves and Koolen pose two further questions about strictly Neumaier graphs which naturally arose from their work:
\begin{enumerate}[(A)]
\item \cite[Question A]{GK18}\label{QA} What is the minimum number of vertices for which there exists a non-strongly-regular, edge-regular graph having a regular clique?
\item \cite[Question B]{GK18}\label{QB} Does there exist a non-strongly-regular, edge-regular graph
having a regular clique with nexus greater than 2?
\end{enumerate}
Indeed, before our work, all known strictly Neumaier graphs had at least $24$ vertices, and contained $m$-regular cliques only for the value $m=1$.
\par
In this paper we answer both of the above questions. First we give some general results on Neumaier graphs and their feasible parameter tuples. In particular, we concentrate on conditions involving parameter tuples that force a Neumaier graph to be strongly regular. We also give a classification of Neumaier graphs with parameters achieving equality in a certain inequality. We then apply these results to determine the smallest strictly Neumaier graph, which turns out to be vertex-transitive and has 16 vertices, valency 9 and a $2$-regular $4$-clique.
\par
Then we present two new infinite sequences of strictly Neumaier graphs. Each of these sequences has first element the unique smallest strictly Neumaier graph. The $i^{\text{th}}$ element of each of these sequences is a strictly Neumaier graph which contains a $2^{i}$-regular clique. In fact, all of these graphs contain a subgraph isomorphic to a clique extension (see the definition of a clique extension in \cite[p.~6]{BCN89}) of the unique smallest strictly Neumaier graph. These constructions show that the nexus of a clique in a strictly Neumaier graph is not bounded above by some constant number. Furthermore, each of the graphs in these sequences has the edge-regular graph parameters of an affine polar graph.

\section{Preliminary}

In this paper we only consider finite, undirected graphs that contain no loops or multiple edges. Let $\Gamma$ be such a graph. We denote by $V(\Gamma)$ the vertex set of $\Gamma$, and $E(\Gamma)$ the edge set of $\Gamma$. For a vertex
$u\in V(\Gamma)$, we define the \emph{neighbourhood }of $u$ in $\Gamma$
to be the set $\Gamma(u)=\{w\in V(\Gamma):uw\in E(\Gamma)\}$. The
\emph{complement }of a graph $\Gamma$, denoted by $\overline{\Gamma}$,
is the graph with vertex set $V(\overline{\Gamma})=V(\Gamma)$ , and for
distinct vertices $u,w\in V(\Gamma)$, we have $uw\in E(\overline{\Gamma})$
if and only of $uw\not\in E(\Gamma)$.

Let $\Gamma$ be a graph and $v=|V(\Gamma)|$. The graph $\Gamma$
is said to be \emph{$k$-regular} if every vertex has neighbourhood
of size $k$. The graph $\Gamma$ is called \emph{regular} if there exists a value $k$ such that $\Gamma$ is $k$-regular. The graph $\Gamma$ is \emph{edge-regular} if it is
non-empty, $k$-regular, and every pair of adjacent vertices have
exactly $\lambda$ common neighbours. Then $\Gamma$ is said to be
edge-regular with \emph{parameters $(v,k,\lambda)$}, and refer to
this as a \emph{parameter tuple}. Denote by $ERG(v,k,\lambda)$ the set of edge-regular graphs with parameters $(v,k,\lambda)$. The graph $\Gamma$ is \emph{co-edge-regular}
if it is non-complete, $k$-regular and every pair of distinct non-adjacent
vertices have exactly $\mu$ common neighbours. Then $\Gamma$ is
said to be co-edge-regular with \emph{parameters} $(v,k,\mu)$. The
graph $\Gamma$ is \emph{strongly regular }if it is both edge-regular
and co-edge-regular. If $\Gamma$ is edge-regular with parameters
$(v,k,\lambda)$, and co-edge-regular with parameters $(v,k,\mu)$,
the graph is said to be strongly regular with \emph{parameters} $(v,k,\lambda,\mu)$.

A \emph{clique} in a graph $\Gamma$ is a set of pairwise adjacent
vertices of $\Gamma$, and a clique of size $s$ is called an $s$\emph{-clique}.
A clique $S$ in a regular graph $\Gamma$ is \emph{regular} if every vertex not in
$S$ is adjacent to the same number $m>0$ of vertices in $S$.\emph{
}In this case we say that $S$ has \emph{nexus} $m$ and is $m$\emph{-regular}. Let us give several examples of strongly regular graphs containing a regular clique.

\medskip
\noindent\textbf{Example 1.}
Let $K_{r\times t}$ be the complete multipartite graph which has $r$ parts of size $t$. Let $S$ be a set consisting of exactly one vertex from each part of $\Gamma$. Then $S$ is a $(r-1)$-regular $r$-clique.

\medskip
\noindent\textbf{Example 2.}
For $n\geq2$, the \emph{square lattice graph} $L_{2}(n)$
has vertex set $\{1,2,...,n\}\times\{1,2,...,n\} $, and two distinct vertices are joined
by an edge precisely when they have the same value at one coordinate.
This graph is strongly regular with parameters $(n^{2},2(n-1),n-2,2)$. Let $S$ be a set consisting of all vertices of $L_{2}(n)$ which have the same fixed value at the same
fixed coordinate. Then $S$ is a $1$-regular $n$-clique.

\medskip
\noindent\textbf{Example 3.}
For $n\geq3$, the \emph{triangular graph} $T(n)$ has vertex
set consisting of the subsets of $\{1,2,...,n\}$ of size $2$, and
two distinct vertices $A,B$ are joined by an edge precisely when
$|A\cap B|=1$. This graph is strongly regular with parameters $({n \choose 2},2(n-2),n-2,4)$. Let $S$ be a set consisting of all
vertices of $T(n)$ which contain a fixed element from $\{1,2,...,n\}$.
Then $S$ is a $2$-regular $(n-1)$-clique.

\medskip
A \emph{Neumaier graph} is a non-complete edge-regular graph which contains a regular clique. We denote by $NG(v,k,\lambda;m,s)$ the set of Neumaier graphs which are edge-regular
with parameters $(v,k,\lambda)$, and contain an $m$-regular $s$-clique, where $s\geq 2$. A \emph{strictly Neumaier graph} is a Neumaier graph which is not strongly regular (the definition of a strictly Neumaier graph is analoguous to the definition of a strictly Deza graph, see \cite{EFHHH99}).

The tuple $(v,k,\lambda)$ is said to be \emph{extremal} if $ERG(v,k,\lambda)$ is non-empty and contains only strongly regular graphs. Similarly, the tuple $(v,k,\lambda;m,s)$ is said to be \emph{extremal} if $NG(v,k,\lambda;m,s)$ is non-empty and contains only strongly regular graphs.

To answer Question \ref{QA}, we collect a series of conditions on the parameters $(v,k,\lambda;m,s)$ that force at least one of the following to occur;\\
{\rm(1)} $ERG(v,k,\lambda)$ is empty.\\
{\rm(2)} $(v,k,\lambda)$ is extremal.  \\
{\rm(3)} $NG(v,k,\lambda;m,s)$ is empty.\\
{\rm(4)} $(v,k,\lambda;m,s)$ is extremal.\\

\medskip
\noindent\textbf{2.1. Edge-regular graphs}
\medskip

First we state simple results concerned with taking the complement
of the graphs we work with.
\begin{lemma}
\label{comps}The following statements hold.\\
{\rm(1)} Suppose $\Gamma$ is a $k$-regular graph. Then $\overline{\Gamma}$ is a $(v-k-1)$-regular graph.\\
{\rm(2)} Suppose $\Gamma$ is an edge-regular graph with parameters $(v,k,\lambda)$.
Then $\overline{\Gamma}$ is co-edge-regular, with parameters $(v,v-k-1,v-2k+\lambda)$.\\
{\rm(3)} Suppose $\Gamma$ is a co-edge-regular graph with parameters $(v,k,\mu)$.
Then $\overline{\Gamma}$ is edge-regular, with parameters $(v,v-k-1,v-2-2k+\mu)$.
\end{lemma}
\begin{corollary}
\label{srgcomp}Let $\Gamma$ be a strongly regular graph with parameters
$(v,k,\lambda,\mu).$ Then $\overline{\Gamma}$ is strongly regular with
parameters $(v,v-k-1,v-2-2k+\mu,v-2k+\lambda)$.
\end{corollary}

The next lemma gives basic properties of an edge-regular graph.

\begin{lemma}[\cite{BCN89}, Chapter 1]
\label{vklam} Let $\Gamma$ be an edge-regular graph with parameters $(v,k,\lambda)$. Then:\\
{\rm(1)} $v>k>\lambda$;\\
{\rm(2)} $v\geq 2k - \lambda$;\\
{\rm(3)} $2$ divides $vk$\label{vk2}; \\
{\rm(4)} $2$ divides $k\lambda$\label{kl2};\\
{\rm(5)} $6$ divides $vk\lambda$\label{vkl6}.\\
\end{lemma}

\noindent\textbf{2.2. Edge-regular graphs with regular cliques}
\medskip

Let $\Gamma$ be an edge-regular graph with parameters $\tau=(v,k,\lambda)$. Our main tool in the investigation of Neumaier graphs is the \emph{clique adjacency polynomial}, which is defined in \cite{S10}, and given by
\[
C_{\tau}(x,y):=x(x+1)(v-y)-2xy(k-y+1)+y(y-1)(\lambda-y+2).
\]
The following theorem uses the clique adjacency polynomial to give a criterion for when any $s$-clique in any graph $\Gamma$ from $ERG(v,k,\lambda)$ is $m$-regular.

\begin{lemma}[\cite{S15}, Theorem 3.1]
\label{cap} Let $\Gamma$ be a graph in $ERG(v,k,\lambda)$ having an $s$-clique $S$, with $s\geq 2$. If $m$ is a positive integer then
\begin{equation}\label{capreg}
C_{\tau}(m-1,s)=C_{\tau}(m,s)=0
\end{equation}
if and only if $S$ is an $m$-regular clique.
\end{lemma}

Further, we list several more tools which we use in the investigation of Neumaier graphs. The next result gives arithmetic conditions on the parameters of a Neumaier graph. By analysing these relations further, we reconstruct $s,m$ as functions of $v,k,\lambda$. The property of these expressions to be integral numbers can then be seen as necessary conditions for an edge-regular graph to contain a regular clique.

\begin{lemma}
\label{count}Let $\Gamma$ be a graph in $NG(v,k,\lambda;m,s)$. Then:\\
{\rm(1)} $(v-s)m=(k-s+1)s$;\\
{\rm(2)} $(k-s+1)(m-1)=(\lambda-s+2)(s-1)$;\\
{\rm(3)} $s$ is the largest root of the polynomial
\begin{equation*}
 (v-2k+\lambda)y^{2} + (k^{2}+3k-\lambda-v(\lambda+2))y + v(\lambda + 1 - k);
\end{equation*}
{\rm(4)} $m$ is the largest root of the polynomial
\begin{equation*}
(v-s)x^{2} - (v-s)x - s(s-1)(\lambda - s + 2).
\end{equation*}
\end{lemma}
\proof
(1) By Theorem \ref{cap}, $C_{\tau}(m,s)=C_{\tau}(m-1,s)=0$ where
$\tau=(v,k,\lambda)$. Then (1) is found by evaluating $0=C_{\tau}(m,s)-C_{\tau}(m-1,s)$.

(2) Substitute (1) into $C_{\tau}(m,s)$.

(3) Multiply expression (2) by $(v-s)$ and use (1) to substitute for $(v-s)m$, we see that $s$ is a root of the polynomial. Note that $v\geq 2k - \lambda$ and $v(\lambda +1-k)\leq 0$ by Lemma \ref{vklam}. This means at most one positive root to the polynomial.

(4) Multiply expression (1) by $(m-1)$ and use (2) to substitute for $(k-s+1)(m-1)$, we see that $m$ is a root of the polynomial. Note that $\lambda-s+2\geq 0$ as an edge in an $s$-clique is in at least $s-2$ triangles of the graph $\Gamma$. This means at most one positive root of the polynomial. $\square$

Now we present a collection of results giving properties of all regular cliques in a Neumaier graph.

\begin{lemma}[\cite{N81}, Theorem 1.1]
Let $\Gamma$ be a graph in $NG(v,k,\lambda;m,s)$. Then:\\
{\rm(1)} the maximum size of a clique in $\Gamma$ is $s$;\\
{\rm(2)} all regular cliques in $\Gamma$ are $m$-regular cliques;\\
{\rm(3)} the regular cliques in $\Gamma$ are precisely the cliques of size $s$.
\end{lemma}

We finish this section by giving a lower bound on the size of a regular clique in a strictly Neumaier graph. We can understand such a result as saying the following: Take a parameter tuple $(v,k,\lambda;m,s)$ where $s$ is less than the bound. Then the parameters $(v,k,\lambda;m,s)$ are extremal.

\begin{lemma}[\cite{GK18}, Proposition 4.2]\label{lam}
Let $\Gamma$ be a strictly Neumaier graph from $NG(v,k,\lambda;m,s)$. Then $s\geq 4$, and consequently, $\lambda\geq2$.
\end{lemma}

\medskip
\noindent\textbf{2.3. Affine polar graphs $VO^+(2e, 2)$}
\medskip

Let $V$ be a $(2e)$-dimensional vector space over a finite field $\mathbb{F}_q$, where $e \ge 2$ and $q$ is a prime power,
provided with the hyperbolic
quadratic form $Q(x) = x_1x_2 + x_3x_4+\ldots+x_{2e-1}x_{2e}$.
The set $Q^+$ of zeroes of $Q$ is called the \emph{hyperbolic quadric}, where $e$ is the maximal dimension of a subspace in $Q^+$.
A \emph{generator} of $Q^+$ is a subspace of maximal dimension $e$ in $Q^+$.

\begin{lemma}[\cite{B16}, Theorem 7.130]\label{PairOfGenerators}
Given $(e-1)$-dimensional subspace $W$ of $Q^+$, there are precisely two generators that contain $W$.
\end{lemma}

Denote by $VO^+(2e,q)$ the graph on $V$ with two vectors $x,y$ being adjacent if and only if $Q(x-y) = 0$. The graph $VO^+(2e,q)$ is known as an \emph{affine polar graph} (see \cite{B16,BH12,BS90}).

\begin{lemma}
A graph $VO^+(2e,q)$ is a vertex transitive strongly regular graph with parameters
\begin{align}
\begin{split}
v &= q^{2e}\\ 
k &= (q^{e-1}+1)(q^e-1)\\
\lambda &= q(q^{e-2}+1)(q^{e-1}-1)+q-2\\
\mu &= q^{e-1}(q^{e-1}+1)
\end{split}
\end{align}
\end{lemma}

Note that $VO^+(2e,q)$ is isomorphic to the graph defined on the set of all $(2\times e)$-matrices over $\mathbb{F}_q$ of the form

\begin{equation}\label{matrices}
\left(
  \begin{array}{cccc}
    x_1 & x_3 & \ldots & x_{2e-1} \\
    x_2 & x_4 & \ldots & x_{2e} \\
  \end{array}
\right)
,
\end{equation}
where two matrices are adjacent if and only if the scalar product of the first and the second rows of their difference is equal to $0$.

\begin{lemma}
There is one-to-one correspondence between cosets of generators of $Q^+$ and maximal cliques in $VO^+(2e,q)$.
\end{lemma}

%\begin{lemma}
%The automorhism group of $VO^+(2e,q)$ acts transitively on the set of maximal cliques in $VO^+(2e,q)$.
%\end{lemma}

\begin{lemma}
Every maximal clique in $VO^+(2e,q)$ is a $q^{e-1}$-regular $q^e$-clique.
\end{lemma}

A \emph{spread} in $VO^+(2e,q)$ is a set of $q^e$ disjoint maximal cliques that correspond to all cosets of a generator.

%\begin{lemma}
%The setwise stabiliser of a given maximal clique in $VO^+(2e,q)$ acts transitively on the set of
%all other maximal cliques from the same spread.
%\end{lemma}

%\begin{lemma}The following statements hold.\\
%{\rm(1)}
%The automorphism group of $VO^+(2e,q)$ acts transitively on the set of all
%$(e-1)$-dimensional isotropic subspaces $W$.\\
%{\rm(2)}
%The automorphism group of $VO^+(2e,q)$ acts transitively on the set of $W$-conjugate generators.\\
%\end{lemma}

\section{Conditions on parameters to be extremal}

We will now give a collection of conditions on parameter tuples to show they are extremal. We first consider the tuples associated with edge-regular graphs, and then consider the tuples associated with Neumaier graphs.

\medskip
\noindent\textbf{3.1. The triple of parameters $(v,k,\lambda)$}
\medskip

When the triple $(v,k,\lambda)$ is extremal, there are no edge-regular graphs in $ERG(v,k,\lambda)$ which are not strongly regular. Thus there is no strictly Neumaier graph with these edge-regular parameters. This fact will be heavily used when analysing the smallest Neumaier graph.

The following Lemma gives a list of sufficient conditions for $(v,k,\lambda)$ to be extremal.

\begin{lemma}
\label{erg} Suppose $ERG(v,k,\lambda)$ is non-empty for some $v,k,\lambda$. Then the triple $(v,k,\lambda)$ is extremal if at least one of the following holds:\\
{\rm(1)} $v=2k-\lambda$\label{miff}.\\
{\rm(2)} $v=2k-\lambda+1$\label{mu1}.\\
{\rm(3)} There is a strongly regular graph with parameters
 $(v,v-k-1,0,v-2k+\lambda)$\label{trifree}.
\end{lemma}
\vspace{-1em}
\proof
(1) These graphs are exactly the graphs $K_{s\times t}$ (see \cite[Theorem 4.1]{S15}).

(2) Take an edge-regular graph with parameters $(v,k,\lambda)$,
with $v-2k+\lambda=1$. By Lemma \ref{comps}, we see that $\overline{\Gamma}$
is co-edge-regular with parameters $(v,v-k-1,1)$. Then by \cite[Lemma 1.1.3]{BCN89},
$\overline{\Gamma}$ is strongly regular. Thus $\Gamma$ is strongly regular.

(3) Let $\Delta$ be a strongly regular graph with parameters\\
 $(v,v-k-1,0,v-2k+\lambda)$. By Corollary \ref{srgcomp}, $\overline{\Delta}$
is strongly regular with parameters $(v,k,\lambda,\mu)$. A standard
counting argument (see \cite[Lemma 1.1.1]{BCN89}) shows
us that $k(k-\lambda-1)=\mu(v-k-1)$.

Now let $u\in V(\Gamma)$. First we partition $V(\Gamma)$ into $V_{1}=\{u\},V_{2}=\Gamma(u)$
and $V_{3}=V(\Gamma)\setminus (\Gamma(u)\cup\{u\})$. Since each vertex
in $V_{2}$ has $k-\lambda-1$ neighbours in $V_{3}$, there are $k(k-\lambda-1)$
edges between $V_{2}$ and $V_{3}$.

Define $b$ as the average number of neighbours a vertex in $V_{3}$
has in $V_{2}$. Then the number of edges between $V_{3}$ and $V_{2}$
is $b(v-k-1).$ Therefore, we have $k(k-\lambda-1)$ is equal to both
$b(v-k-1)$ and $\mu(v-k-1)$, so $b=\mu$.

Let $w\in V_{3}$. The number of neighbours of $w$ in $V_{2}$ is
at least $k-(v-k-2)=\mu$, as $|V_{3}|=v-k-1$ . As $b=\mu$ is the
average of numbers at least as big as $b$, they must all equal $b$.
This means the number of common neighbours of $u$ and $w$ in $\Gamma$
is exactly $\mu$, and so $\Gamma$ is strongly regular. $\square$

\medskip
\noindent\textbf{3.2. The quintuple of parameters $(v,k,\lambda;m,s)$}
\medskip

Next we will give a necessary condition for the existence of a graph in $NG(v,k,\lambda;m,s)$, in the form of an inequality that is linear in the parameters $k,\lambda,m,s$. When equality is achieved, we show that the parameters $(v,k,\lambda;m,s)$ are extremal.

Firstly, we give a useful Lemma involving Neumaier graphs where the neighbourhood of any vertex has a certain structure.

\begin{lemma}\label{genregvtx}
Let $\Gamma$ be a graph from $NG(v,k,\lambda;m,s)$. Further suppose that every vertex
in $\Gamma$ has neighbourhood consisting of $l$ vertex disjoint
cliques of size $s-1$. Then $\Gamma$ is strongly regular, with parameters
$v=s+(l-1)(s-1)s/m$, $k=l(s-1)$, $\lambda=(m-1)(l-1)+s-2$ and $\mu=lm$.
\end{lemma}
\proof
Take any vertex $u\in V(\Gamma)$ and $w\not\in\Gamma(u)$. The neighbourhood
of $u$ consists of disjoint $(s-1)$-cliques. Together with $u$
each of these cliques define an $s$-clique. These cliques are necessarily
$m$-regular by Theorem \ref{cap}. Thus $w$ is adjacent to $m$ vertices
in each of these cliques, and has exactly $lm$ neighbours in common
with $u$. This proves $\Gamma$ is strongly regular with $\mu=lm$.

The formulae for $k$ and $\lambda$ can be derived by simple counting
arguments. Then for $v$, we use Proposition \ref{count}. $\square$

Now we give the inequality of the parameters $(v,k,\lambda)$ of a Neumaier graph. In the equality case, we show we are in a situation covered by the Lemma \ref{genregvtx}, proving that $(v,k,\lambda;m,s)$ is extremal.

\begin{theorem}
\label{ksm}Let $\Gamma$ be a graph from $NG(v,k,\lambda;m,s)$. Then
\[
k-\lambda-s+m-1\geq0\tag*{(*)}
\]
 Equality holds if and only if every vertex in $\Gamma$ has a neighbourhood
consisting of two vertex disjoint $(s-1)$-cliques. In this case,
$\Gamma$ is a complete graph or strongly regular with $v=s+(s(s-1)/m)$, $k=2(s-1)$,
$\lambda=s+m-3$ and $\mu=2m$.
\end{theorem}

\proof
Let $u\in S$, and consider $w\in V(\Gamma)\setminus S$, with $uw\in E(\Gamma)$.
We know that $u$ has $k-s$ other neighbours in $V(\Gamma)\setminus S$,
and $w$ has $m-1$ neighbours in $S\setminus\{u\}$. Thus $u$ and
$w$ have exactly $m-1$ common neighbours in $S$, and at most $k-s$
common neighbours in $V(\Gamma)\setminus S$. As $u,w$ have exactly
$\lambda$ common neighbours, we must have $\lambda\leq k-s+m-1$.

When equality holds, we see that $w$ must be adjacent to all neighbours
of $u$ in $V(\Gamma)\setminus S$. By repeating the argument for
all other edges $uz$, with $z\in V(\Gamma)\setminus S$, we see
that $u$ has a neighbourhood consisting of two vertex disjoint cliques.

By Proposition \ref{count} (2) and $k=\lambda+s-m+1$, we deduce
that $(\lambda-s-m+3)(s-m)=0$. If $s=m$, $\Gamma$ is necessarily
complete. Otherwise, $\lambda=s+m-3$ and $k=2(s-1)$. This proves
that for all $u\in S$, $u$ has a neighbourhood consisting of two
vertex disjoint $(s-1)$-cliques.

Now take a vertex $u\in V(\Gamma)\setminus S$. As $m\geq1$, $u$
is adjacent to a vertex $w\in S$. As the neighbourhood of $w$ consists
of $(s-1)$-cliques, $u$ is contained in $S^{'}$, which is one of
these $(s-1)$-cliques. Then $S=S^{'}\cup\{w\}$ is an $s$-clique
that contains $u$. Thus, we have proven every vertex is contained
in an $s$-clique.

By Theorem \ref{cap}, any $s$-clique is necessarily $m$-regular. So we
can apply the above argument to show that any vertex $u\in V(\Gamma)$
has a neighbourhood consisting of $2$ vertex disjoint $(s-1)$-cliques
in $\Gamma$. The result then follows from Lemma \ref{genregvtx}. $\square$

\medskip
\noindent\textbf{3.2.1. Classifying the graphs in the equality case}
\medskip

In fact, we can give a full description of all Neumaier graphs with parameters satisfying equality in (*).

\begin{theorem}
Let $\Gamma$ be a graph from $NG(v,k,\lambda;m,s)$, where $k-\lambda-s+m-1=0$. Then $\Gamma$ is one of the following strongly regular graphs:\\
{\rm(1)} the square lattice graph $L_{2}(s)$;\\
{\rm(2)} the triangular graph $T(s+1)$, where $s\geq3$;\\
{\rm(3)} the complete $s$-partite graph $K_{s\times2}$, with parts of size $2$.
\end{theorem}

We prove this theorem by taking cases on the value of $m$. We start with the case $m=1$.

\begin{lemma}
\label{latgraph}
Let $\Gamma$ be a graph from $NG(v,k,\lambda;m,s)$, where $k-\lambda-s+m-1=0$ and $m=1$. Then $\Gamma$ is isomorphic to the 
square lattice graph $L_{2}(s)$.
\end{lemma}
\proof
By Theorem \ref{ksm}, $\Gamma$ is strongly regular with parameters $(s^{2},2(s-1),s-2,2)$. Any strongly regular graph with parameters $(s^{2},2(s-1),s-2,2)$ must be isomorphic to $L_{2}(s)$,
unless $s=4$ (see \cite{S59b}). In this case, there is only one strongly regular graph that is not isomorphic to $L_{2}(4)$, called the \emph{Shrikhande graph}. This graph does not contain a regular clique (see \cite{C59}). $\square$

\medskip
Next we consider the case $m=2$.

\begin{lemma}
\label{trigraph}Let $\Gamma$ be a graph from $NG(v,k,\lambda;m,s)$, where $k-\lambda-s+m-1=0$ and $m=2$. Then $\Gamma$ is isomorphic to the triangular graph $T(s+1)$.
\end{lemma}
\proof
By Theorem \ref{ksm}, $\Gamma$ is strongly regular with parameters $(s^{2},2(s-1),s-2,2)$. Any strongly regular graph with parameters $({s+1 \choose 2},2(s-1),s-1,4)$ must be isomorphic to $T(s+1)$, unless $s=7$ (see \cite{C58,H60,S59a} or \cite{C59}). In this case, there are only three strongly regular graphs that are not isomorphic to $T(8)$, called the \emph{Chang graphs}. Each of these do not contain a regular clique (see \cite{C59}). $\square$

\medskip
Now we only need to consider the case $m\geq3$. For this case, we can show that $m$ is particularly large with respect to $s$, which forces the graph to be isomorphic to $K_{s\times 2}$.

\begin{lemma}
\label{matgr} Let $\Gamma$ be a graph from $NG(v,k,\lambda;m,s)$, where $k-\lambda-s+m-1=0$ and $m\geq3$. Then $\Gamma$ is isomorphic to $K_{s\times2}$.
\end{lemma}

\proof

We will first show that $m\geq1+s/2$.

Let $\Gamma$ be a graph in $NG(v,k,\lambda)$. Take a subset $S\subset V(\Gamma)$, $S=\{u_{1},u_{2},\ldots,u_{s}\}$, where $S$ a $m$-regular $s$-clique in $\Gamma$. Without loss of generality, let $w\in V(\Gamma)\setminus S$, with $\{u_{1},u_{2},u_{3}\}\subseteq\Gamma(w)\cap S$.
Note that by the equality case of Theorem \ref{ksm}, $w$ is adjacent to
all neighbours of $u_{1},u_{2},u_{3}$ in $V(\Gamma)\setminus S$.

As $\Gamma$ is $k$-regular, we have $|N(u_{i})\cap V(\Gamma)\setminus S|=k-s+1$
for all $i\in\{1,\ldots,s\}$. Also we must have $|N(u_{i})\cap N(u_{j})\cap V(\Gamma)\setminus S|=\lambda-s+2$
for all $i,j\in\{1,\ldots,s\}$. Thus we have
\begin{equation*}
\begin{split}
|(N(u_{1})\cup N(u_{2})\cup N(u_{3}))\cap V(\Gamma)\setminus S|&=3(k-s+1)-3(\lambda-s+2)\\
                 &+|(N(u_{1})\cap N(u_{2})\cap N(u_{3}))\cap V(\Gamma)\setminus S|
\end{split}
\end{equation*}
We see that $w$ is adjacent to at least $m+3(k-\lambda-1)$ vertices.
Therefore, this has to be less than $k$. By using $\lambda=k-s+m-1$
and $k=2(s-1)$ (by Theorem \ref{ksm}), we get that $m\geq1+s/2$.

Let $u,w\in V(\Gamma)\setminus S$. As $S$ is $m$-regular and $m>s/2$,
there must exist a $i\in S$ such that $iu,iw\in E(\Gamma)$. We also
know that the neighbourhood of $i$ in $V(\Gamma)\setminus S$ is
a clique, so $uw\in E(\Gamma)$. Thus we have proven $V(\Gamma)\setminus S$
is a clique in $\Gamma$.

By maximality of $S$, we must have $|V(\Gamma)\setminus S|\leq s$.
Also, because $k=2(s-1)$, we must have $|V(\Gamma)\setminus S|\geq s-1$.
As $\Gamma$ is non-complete, we have $|V(\Gamma)\setminus S|=s$,
and $m=s-1$. By Theorem \ref{ksm}, $v=2s,\lambda=2(s-2)$. Applying Proposition \ref{erg},
we get the result. $\square$

\section{Small extremal parameter tuples and the smallest strictly Neumaier graph}

The following tables list all tuples $v,k,\lambda,s,m$ of integers,
such that the following hold:
\begin{enumerate}
\item $0<k<v-1,\ v\leq24,\ 0\leq\lambda<k,\ 2\leq s\leq\lambda+2$ and $m\geq1$.
\item $2$ divides both $vk$ and $k\lambda$, and $6$ divides $vk\lambda$
(see Lemma (\ref{vklam})).
\item $C_{\tau}(m-1,s)=C_{\tau}(m,s)=0$, where $\tau=(v,k,\lambda)$ (see
Lemma \ref{cap}).
\end{enumerate}
These tables were obtained by a straightforward computation using
GAP \cite{T17}. All calculations were exact and
took a total of about 20 CPU milliseconds on a desktop PC.

Thus, if there is a Neumaier graph from $NG(v,k,\lambda;m,s)$ such that $v\leq24$,
then the tuple $v,k,\lambda,m,s$ appears in our tables. The last
columns of our tables display a result which proves that the tuple $(v,k,\lambda)$ or the tuple $(v,k,\lambda;m,s)$ is extremal, or the symbol `-' otherwise. For example, L\ref{erg} (1) refers to Lemma \ref{erg} part (1), and T\ref{ksm} refers to Theorem \ref{ksm}.

\begin{table}[H]
\begin{minipage}[t]{0.4\columnwidth}%
\begin{tabular}{c|c|c|c|c|c}
$v$ & $k$ & $\lambda$ & $m$ & $s$ & result\tabularnewline
\hline
4 & 2 & 0 & 1 & 2 & L\ref{erg} (1)\tabularnewline
\hline
6 & 3 & 0 & 1 & 2 & L\ref{erg} (1)\tabularnewline
\cline{2-6}
 & 4 & 2 & 2 & 3 & L\ref{erg} (1)\tabularnewline
\hline
8 & 4 & 0 & 1 & 2 & L\ref{erg} (1)\tabularnewline
\cline{2-6}
 & 6 & 4 & 3 & 4 & L\ref{erg} (1)\tabularnewline
\hline
9 & 4 & 1 & 1 & 3 & L\ref{lam}\tabularnewline
\cline{2-6}
 & 6 & 3 & 2 & 3 & L\ref{erg} (1)\tabularnewline
\hline
10 & 5 & 0 & 1 & 2 & L\ref{erg} (1)\tabularnewline
\cline{2-6}
 & 6 & 3 & 2 & 4 & L\ref{erg} (3)\tabularnewline
\cline{2-6}
 & 8 & 6 & 4 & 5 & L\ref{erg} (1)\tabularnewline
\hline
12 & 5 & 2 & 1 & 4 & T\ref{ksm}\tabularnewline
\cline{2-6}
 & 6 & 0 & 1 & 2 & L\ref{erg} (1)\tabularnewline
\cline{3-6}
 &  & 4 & 1 & 6 & T\ref{ksm}\tabularnewline
\cline{2-6}
 & 8 & 4 & 2 & 3 & L\ref{erg} (1)\tabularnewline
\cline{2-6}
 & 9 & 6 & 3 & 4 & L\ref{erg} (1)\tabularnewline
\cline{2-6}
 & 10 & 8 & 5 & 6 & L\ref{erg} (1)\tabularnewline
\end{tabular}%
\end{minipage}\hfill{}%
\begin{minipage}[t]{0.4\columnwidth}%
\begin{tabular}{c|c|c|c|c|c}
$v$ & $k$ & $\lambda$ & $m$ & $s$ & result\tabularnewline
\hline
14 & 7 & 0 & 1 & 2 & L\ref{erg} (1)\tabularnewline
\cline{2-6}
 & 9 & 6 & 3 & 7 & T\ref{ksm}\tabularnewline
\cline{2-6}
 & 12 & 10 & 6 & 7 & L\ref{erg} (1)\tabularnewline
\hline
15 & 6 & 1 & 1 & 3 & L\ref{lam}\tabularnewline
\cline{3-6}
 &  & 3 & 1 & 5 & T\ref{ksm}\tabularnewline
\cline{2-6}
 & 8 & 4 & 2 & 5 & T\ref{ksm}\tabularnewline
\cline{2-6}
 & 10 & 5 & 2 & 3 & L\ref{erg} (1)\tabularnewline
\cline{3-6}
 &  & 6 & 3 & 5 & L\ref{mu1}\tabularnewline
\cline{2-6}
 & 12 & 9 & 4 & 5 & L\ref{erg} (1)\tabularnewline
\hline
16 & 6 & 2 & 1 & 4 & T\ref{ksm}\tabularnewline
\cline{2-6}
 & 8 & 0 & 1 & 2 & L\ref{erg} (1)\tabularnewline
\cline{3-6}
 &  & 6 & 1 & 8 & T\ref{ksm}\tabularnewline
\cline{2-6}
 & 9 & 4 & 2 & 4 & -\tabularnewline
\cline{2-6}
 & 10 & 6 & 3 & 6 & L\ref{erg} (3)\tabularnewline
\cline{2-6}
 & 12 & 8 & 3 & 4 & L\ref{erg} (1)\tabularnewline
\cline{2-6}
 & 14 & 12 & 7 & 8 & L\ref{erg} (1)\tabularnewline
\multicolumn{1}{c}{} & \multicolumn{1}{c}{} & \multicolumn{1}{c}{} & \multicolumn{1}{c}{} & \multicolumn{1}{c}{} & \tabularnewline
\end{tabular}%
\end{minipage}

\caption*{Table 1. : Possible parameters of Neumaier graphs on $v\leq16$ vertices }
\end{table}

We see that Table 1 rules out all possible parameter
tuples $(v,k,\lambda;m,s)$ for a strictly Neumaier graph when $v<16$. Further, the table shows
that any strictly Neumaier graph on $16$ vertices is from  $NG(16,9,4;2,4)$. A graph in $NG(16,9,4;2,4)$ is given in the next section.
So Table 1 and this graph give the answer to Question
\ref{QA}.

\begin{table}[H]
\begin{minipage}[t]{0.4\columnwidth}%
\begin{tabular}{c|c|c|c|c|c}
$v$ & $k$ & $\lambda$ & $m$ & $s$ & result\tabularnewline
\hline
18 & 7 & 4 & 1 & 6 & T\ref{ksm}\tabularnewline
\cline{2-6}
 & 9 & 0 & 1 & 2 & L\ref{erg} (1)\tabularnewline
\cline{2-6}
 & 12 & 6 & 2 & 3 & L\ref{erg} (1)\tabularnewline
\cline{2-6}
 & 15 & 12 & 5 & 6 & L\ref{erg} (1)\tabularnewline
\cline{2-6}
 & 16 & 14 & 8 & 9 & L\ref{erg} (1)\tabularnewline
\hline
20 & 10 & 0 & 1 & 2 & L\ref{erg} (1)\tabularnewline
\cline{2-6}
 & 15 & 10 & 3 & 4 & L\ref{erg} (1)\tabularnewline
\cline{2-6}
 & 16 & 12 & 4 & 5 & L\ref{erg} (1)\tabularnewline
\cline{2-6}
 & 18 & 16 & 9 & 10 & L\ref{erg} (1)\tabularnewline
\hline
21 & 8 & 1 & 1 & 3 & L\ref{lam}\tabularnewline
\cline{3-6}
 &  & 5 & 1 & 7 & T\ref{ksm}\tabularnewline
\cline{2-6}
 & 10 & 5 & 2 & 6 & T\ref{ksm}\tabularnewline
\cline{2-6}
 & 12 & 7 & 3 & 7 & T\ref{ksm}\tabularnewline
\cline{3-6}
 &  & 8 & 3 & 9 & T\ref{ksm}\tabularnewline
\cline{2-6}
 & 14 & 7 & 2 & 3 & L\ref{erg} (1)\tabularnewline
\cline{3-6}
 &  & 9 & 4 & 7 & -\tabularnewline
\cline{2-6}
 & 15 & 10 & 4 & 6 & L\ref{erg} (2)\tabularnewline
\cline{2-6}
 & 16 & 12 & 6 & 9 & L\ref{erg} (2)\tabularnewline
\cline{2-6}
 & 18 & 15 & 6 & 7 & L\ref{erg} (1)\tabularnewline
\end{tabular}%
\end{minipage}\hfill{}%
\begin{minipage}[t]{0.4\columnwidth}%
\begin{tabular}{c|c|c|c|c|c}
$v$ & $k$ & $\lambda$ & $m$ & $s$ & result\tabularnewline
\hline
22 & 11 & 0 & 1 & 2 & L\ref{erg} (1)\tabularnewline
\cline{2-6}
 & 12 & 5 & 2 & 4 & -\tabularnewline
\cline{2-6}
 & 14 & 9 & 4 & 8 & T\ref{ksm}\tabularnewline
\cline{2-6}
 & 16 & 12 & 6 & 11 & T\ref{ksm}\tabularnewline
\cline{2-6}
 & 20 & 18 & 10 & 11 & L\ref{erg} (1)\tabularnewline
\hline
24 & 8 & 2 & 1 & 4 & -\tabularnewline
\cline{3-6}
 &  & 4 & 1 & 6 & T\ref{ksm}\tabularnewline
\cline{2-6}
 & 9 & 6 & 1 & 8 & T\ref{ksm}\tabularnewline
\cline{2-6}
 & 12 & 0 & 1 & 2 & L\ref{erg} (1)\tabularnewline
\cline{3-6}
 &  & 10 & 1 & 12 & T\ref{ksm}\tabularnewline
\cline{2-6}
 & 16 & 8 & 2 & 3 & L\ref{erg} (1)\tabularnewline
\cline{2-6}
 & 18 & 12 & 3 & 4 & L\ref{erg} (1)\tabularnewline
\cline{3-6}
 &  & 15 & 6 & 16 & T\ref{ksm}\tabularnewline
\cline{2-6}
 & 20 & 16 & 5 & 6 & L\ref{erg} (1)\tabularnewline
\cline{3-6}
 &  & 17 & 10 & 16 & T\ref{ksm}\tabularnewline
\cline{2-6}
 & 21 & 18 & 7 & 8 & L\ref{erg} (1)\tabularnewline
\cline{2-6}
 & 22 & 20 & 11 & 12 & L\ref{erg} (1)\tabularnewline
\multicolumn{1}{c}{} & \multicolumn{1}{c}{} & \multicolumn{1}{c}{} & \multicolumn{1}{c}{} & \multicolumn{1}{c}{} & \tabularnewline
\multicolumn{1}{c}{} & \multicolumn{1}{c}{} & \multicolumn{1}{c}{} & \multicolumn{1}{c}{} & \multicolumn{1}{c}{} & \tabularnewline
\end{tabular}%
\end{minipage}

\caption*{Table 2. : Possible parameters of Neumaier graphs on
 $16<v\leq24$ vertices }
\end{table}

Table 1 and 2 together show that $(24,8,2)$
is the only possible parameter tuple for a strictly Neumaier graph containing
a $1$-regular clique when $v\leq24$.
%We can then rule out the existence of strictly Neumaier graphs with parameter tuples $(21,14,9;4,7)$ and $(22,12,5;2,4)$ computationally.

Finally, direct computations can show that there is only strictly Neumaier graph in $NG(16,9,4;2,4)$, up to isomorphism,
and there are no strictly Neumaier graphs with parameter tuples $(21,14,9;4,7)$ and $(22,12,5;2,4)$.
Let us explain some ideas of the computations. We fix a subgraph induced by vertices of a clique with given size.
Then we exhaust all regular graphs such that the fixed clique is regular with given nexus.
Using MAGMA \cite{BCP97}, we find that the graphs in $NG(16,9,4;2,4)$ are isomorphic pairwise.

Thus we have found that any strictly Neumaier graph on at most $24$ vertices must have parameters $(16,9,4;2,4)$ or $(24,8,2;1,4)$.

\medskip
\noindent\textbf{4.1. Vertex-transitive strictly Neumaier graphs}
\medskip

The authors discovered the smallest strictly Neumaier graph independently, using completely different approaches.

Goryainov and Panasenko
were looking for strictly Neumaier graphs that admit a partition into regular cliques and used this pattern
for computer searching.

Evans found the graph in a collection of vertex-transitive edge-regular graphs received from Gordon Royle. Holt and Royle have recently enumerated all transitive permutation groups of degree at most $47$ \cite{HR17}. From this, Royle was able to enumerate all vertex-transitive edge-regular graphs on less than $47$ vertices.

Thus we also find all vertex-transitive strictly Neumaier graphs on at most $47$ vertices using the enumeration \cite{HR17}. We list the parameters of all vertex-transitive strictly Neumaier graph on at most $47$ vertices, and the number of vertex-transitive strictly Neumaier graphs with these parameters.\\
{\rm{(1)} $1$ graph with parameters $(16,9,4;2,4)$.\\
{\rm{(2)} $4$ graphs with parameters $(24,8,2;1,4)$.\\
{\rm{(3)} $2$ graphs with parameters $(28,9,2;1,4)$.\\
{\rm{(4)} $1$ graph with parameters $(40,12,2;1,4)$.\\
We note that the four vertex-transitive strictly Neumaier graphs in $NG(24,8,2;1,4)$ appear in \cite{GS14}. They come about in a search for Deza
graphs, which are a certain generalisation of strongly regular graphs.

\medskip
\section{Two constructions and two generalisations of the smallest strictly Neumaier graph}

In this section, we will construct two sequences of strictly Neumaier graphs that generalise the smallest strictly Neumaier graph. The motivation behind both constructions is as follows.

Consider a graph $\Gamma$, and two disjoint subsets $S,T$ of the vertices of $\Gamma$. Now we introduce an important operation on the graph $\Gamma$. For all vertices in $v$ in $S$, we do the following. First take $N=N(v)\cap T$, the neighbours of $v$ in $T$, and $M=T\setminus N$. Then delete all edges $vu$ where $u$ is in $N$, and insert all edges $vu$ where $u$ is in $M$. We will call this operation a \emph{switching} of the edges between $S$ and $T$ in the graph $\Gamma$.

Note that the smallest Neumaier graph contains disjoint $2$-regular $4$-cliques. A switching between any distinct pair of these cliques will not change the fact that they are $2$-regular. Therefore, if we could find a strongly regular graph with these parameters and containing disjoint $2$-regular $4$-cliques, we could hope that the smallest Neumaier graph is the result of switching edges between them.

In the following subsections we will see that the smallest Neumaier graph is the result of two consecutive switchings of the graph $VO^+(4,2)$. We then generalise our switchings to the graphs $VO^+(2e,2)$ for larger $e$, and construct infinite sequences of strictly Neumaier graphs with the same edge-regular parameters as $VO^+(2e,2)$. From now on, we will denote $VO^+(2e,2)$ as the graph $\Gamma_{e}$. Throughout this section we use matrix notation with stars `*' as entries, which denotes the set of corresponding matrices
where the stars take all possible values independently.

\medskip
\noindent\textbf{5.1. The first construction of the smallest Neumaier graph}
\medskip

Consider the $1$-dimensional subspace
$$
W =
\left(
  \begin{array}{cccc}
    \ast & 0 \\
    0  & 0 \\
  \end{array}
\right).
$$
According to Lemma \ref{PairOfGenerators}, the subspace $W$ is contained in the two generators
$$
W_1
=
\left(
  \begin{array}{cccc}
    \ast & \ast \\
    0  & 0 \\
  \end{array}
\right)
\text{ and }
W_2
=
\left(
  \begin{array}{cccc}
    \ast & 0\\
    0  & \ast  \\
  \end{array}
\right).
$$
Take the vector
$$
v =
\left(
  \begin{array}{cccc}
    0 & 0\\
    1  & 0  \\
  \end{array}
\right)
$$
and consider the cosets
$$
v + W_1
=
\left(
  \begin{array}{cccc}
    \ast & \ast \\
    1  & 0 \\
  \end{array}
\right),
$$
$$
v + W_2
=
\left(
  \begin{array}{cccc}
    \ast & 0\\
    1  & \ast  \\
  \end{array}
\right),
$$
whose intersection is
$$
v + W
=
\left(
  \begin{array}{cccc}
    \ast & 0\\
    1  & 0  \\
  \end{array}
\right).
$$
In this setting, the adjacency matrix of the affine polar graph $\Gamma_{2}=VO^+(4,2)$ can be seen in Figure 1. The graph $\Gamma_{2}$ is isomorphic to the complement of the square lattice graph $L_2(4)$.

\begin{table}[H]

$$
  \begin{array}{c|@{}c@{}@{}c@{}@{}c@{}@{}c@{}@{}c@{}@{}c@{}|@{}c@{}@{}c@{}@{}c@{}@{}c@{}@{}c@{}@{}c@{}|@{}c@{}@{}c@{}@{}c@{}@{}c@{}}
  & ~
  \begin{array}{@{}c@{}@{}c@{}}
    01 \\
    00 \\
  \end{array} ~
& ~
  \begin{array}{@{}c@{}@{}c@{}}
    11 \\
    00 \\
  \end{array} ~

 &~
   \begin{array}{@{}c@{}@{}c@{}}
    00 \\
    00 \\
  \end{array}~
  & ~
    \begin{array}{@{}c@{}@{}c@{}}
    10 \\
    00 \\
  \end{array} ~
   & ~
     \begin{array}{@{}c@{}@{}c@{}}
    00 \\
    01 \\
  \end{array} ~
    & ~
      \begin{array}{@{}c@{}@{}c@{}}
    10 \\
    01 \\
  \end{array} ~~
     & ~
  \begin{array}{@{}c@{}@{}c@{}}
    01 \\
    10 \\
  \end{array} ~
      & ~
          \begin{array}{@{}c@{}@{}c@{}}
    11 \\
    01 \\
  \end{array} ~
       & ~
            \begin{array}{@{}c@{}@{}c@{}}
    00 \\
    10 \\
  \end{array} ~
        & ~
          \begin{array}{@{}c@{}@{}c@{}}
    10 \\
    10 \\
  \end{array} ~
         & ~
   \begin{array}{@{}c@{}@{}c@{}}
    00 \\
    11 \\
  \end{array} ~
          & ~
         \begin{array}{@{}c@{}@{}c@{}}
    10 \\
    11 \\
  \end{array} ~~
           & ~
     \begin{array}{@{}c@{}@{}c@{}}
    01 \\
    01 \\
  \end{array} ~
        & ~
            \begin{array}{@{}c@{}@{}c@{}}
    11 \\
    01 \\
  \end{array} ~
         & ~
             \begin{array}{@{}c@{}@{}c@{}}
    01 \\
    11 \\
  \end{array} ~
          & ~
              \begin{array}{@{}c@{}@{}c@{}}
    11 \\
    11 \\
  \end{array}  ~   \\
  \hline

    \begin{array}{@{}c@{}@{}c@{}}
    01 \\
    00 \\
  \end{array}
& ~0 & ~1 & ~1 &   ~1 &  ~0 &  0 & ~\color{red}{1} & ~\color{red}{0} & ~\color{red}{1} & ~\color{red}{0} &  ~0 &  1 &  ~1 &  ~1 &  ~1 &  ~0         \\
  \begin{array}{@{}c@{}@{}c@{}}
    11 \\
    00 \\
  \end{array}

 &   ~1  &  ~0  &  ~1  &  ~1  &  ~0  &  0  & ~\color{red}{0}  & ~\color{red}{1}  & ~\color{red}{0}  & ~\color{red}{1} &  ~1  &  0  &   ~1 &   ~1  &  ~0   & ~1               \\
   \begin{array}{@{}c@{}@{}c@{}}
    00 \\
    00 \\
  \end{array}
  &    ~1  &  ~1  &  ~0  &  ~1  &  ~1  &  1  &  ~\color{red}{1}  &  ~\color{red}{0} &  ~1  &  ~0 &  ~\color{red}{1}  &  \color{red}{0}  &   ~0  &  ~0  &  ~0  &  ~1             \\
    \begin{array}{@{}c@{}@{}c@{}}
    10 \\
    00 \\
  \end{array}
   &    ~1  &  ~1  &  ~1  &  ~0  &  ~1  &  1  &  ~\color{red}{0}  & ~\color{red}{1} & ~0  &  ~1  &  ~\color{red}{0}  &  \color{red}{1}  &   ~0  &  ~0  &  ~1  &  ~0            \\
     \begin{array}{@{}c@{}@{}c@{}}
    00 \\
    01 \\
  \end{array}
    &  ~0  &  ~0  &  ~1  &  ~1  &  ~0  &  1  &  ~0  &  ~1  & ~\color{red}{1}  &  ~\color{red}{0}  &  ~\color{red}{1}  &  \color{red}{0}  &   ~1  &  ~1  &  ~1  &  ~0             \\
      \begin{array}{@{}c@{}@{}c@{}}
    10 \\
    01 \\
  \end{array}
     &  ~0  &  ~0  &  ~1  &  ~1  &  ~1  &  0  &  ~1  &  ~0  &  ~\color{red}{0}  &  ~\color{red}{1}  &  ~\color{red}{0}  &  \color{red}{1}  &  ~1  &  ~1  &  ~0  &  ~1             \\
     \hline
  \begin{array}{@{}c@{}@{}c@{}}
    01 \\
    10 \\
  \end{array}
      &   ~\color{red}{1}  &  ~\color{red}{0}  &  ~\color{red}{1}  &  ~\color{red}{0}  &  ~0  &  1  &  ~0  &  ~1  &  ~1  &  ~1  &  ~0  &  ~0  &  ~1  &  ~0  &  ~1  &  ~1           \\
          \begin{array}{@{}c@{}@{}c@{}}
    11 \\
    01 \\
  \end{array}
       &  ~\color{red}{0}  &  ~\color{red}{1}  &  ~\color{red}{0}  &  ~\color{red}{1}  &  ~1  &  0  &  ~1  &  ~0  &  ~1  &  ~1  &  ~0  &  ~0  &  ~0  &  ~1  &  ~1  &  ~1           \\
            \begin{array}{@{}c@{}@{}c@{}}
    00 \\
    10 \\
  \end{array}
        &   ~\color{red}{1}  &  ~\color{red}{0}  &  ~1  &  ~0  &  ~\color{red}{1}  &  \color{red}{0}  &  ~1  &  ~1  &  ~0  &  ~1  &  ~1  &  ~1  &  ~0  &  ~1  &  ~0  &  ~0             \\
          \begin{array}{@{}c@{}@{}c@{}}
    10 \\
    10 \\
  \end{array}
         &  ~\color{red}{0}  &  ~\color{red}{1}  &  ~0  &  ~1  &  ~\color{red}{0}  &  \color{red}{1}  &  ~1  &  ~1  &  ~1  &  ~0  &  ~1  &  ~1  &  ~1  &  ~0  &  ~0  &  ~0                \\
   \begin{array}{@{}c@{}@{}c@{}}
    00 \\
    11 \\
  \end{array}
          &  ~0  &  ~1  &  ~\color{red}{1}  &  ~\color{red}{0}  &  ~\color{red}{1}  &  \color{red}{0}  &  ~0  &  ~0  &  ~1  &  ~1  &  ~0  &  ~1  &  ~1  &  ~0  &  ~1  &  ~1              \\
         \begin{array}{@{}c@{}@{}c@{}}
    10 \\
    11 \\
  \end{array}
           & ~1  &  ~0  &  ~\color{red}{0}  &  ~\color{red}{1}  &  ~\color{red}{0}  &  \color{red}{1}  &  ~0  &  ~0  &  ~1  &  ~1  &  ~1  &  ~0  &  ~0  &  ~1  &  ~1  &  ~1           \\
     \hline
     \begin{array}{@{}c@{}@{}c@{}}
    01 \\
    01 \\
  \end{array}
        &  ~1  &  ~1  &  ~0  &  ~0  &  ~1  &  1  &  ~1  &  ~0  &  ~0  &  ~1  &  ~1  &  ~0  &   ~0  &  ~1  &  ~1  &  ~0                 \\
            \begin{array}{@{}c@{}@{}c@{}}
    11 \\
    01 \\
  \end{array}
         &  ~1  &  ~1  &  ~0  &  ~0  &  ~1  &  1  &  ~0  &  ~1  &  ~1  &  ~0  &  ~0  &  ~1  &  ~1  &  ~0  &  ~0  &  ~1                \\
             \begin{array}{@{}c@{}@{}c@{}}
    01 \\
    11 \\
  \end{array}
          &   ~1  &  ~0  &  ~0  &  ~1  &  ~1  &  0  &  ~1  &  ~1  &  ~0  &  ~0  &  ~1  &  ~1  &  ~1  &  ~0  &  ~0  &  ~1               \\
              \begin{array}{@{}c@{}@{}c@{}}
    11 \\
    11 \\
  \end{array}
  & ~0  &  ~1  &  ~1  &  ~0  &  ~0  &  1  &  ~1  &  ~1  &  ~0  &  ~0  &  ~1  &  ~1  &  ~0  &  ~1  &  ~1  &  ~0    \\
  \end{array}
$$

\caption*{Figure 1: The adjacency matrix, $A_{2}$, of $\Gamma _{2}=VO^{+}(4,2)$ }
\end{table}
We note that switching edges between the cliques $W_1$, $v+W_1$ gives a graph isomorphic to the complement of the Shrikhande
graph. The switching of edges between the cliques $W_1$, $v+W_1$ and then between the cliques $W_2$, $v+W_2$
is equivalent to inverting the red entries. This gives the strictly Neumaier graph $\Gamma_{2,1}$ on 16 vertices,
whose adjacency matrix is presented in Figure 2.

\begin{table}[H]
$$
  \begin{array}{@{}c@{}|@{}c@{}@{}c@{}@{}c@{}@{}c@{}@{}c@{}@{}c@{}|@{}c@{}@{}c@{}@{}c@{}@{}c@{}@{}c@{}@{}c@{}|@{}c@{}@{}c@{}@{}c@{}@{}c@{}|@{}c@{}}
  & ~
  \begin{array}{@{}c@{}@{}c@{}}
    01 \\
    00 \\
  \end{array} ~
&  ~
  \begin{array}{@{}c@{}@{}c@{}}
    11 \\
    00 \\
  \end{array}~

 &~
   \begin{array}{@{}c@{}@{}c@{}}
    00 \\
    00 \\
  \end{array}~
  &~
    \begin{array}{@{}c@{}@{}c@{}}
    10 \\
    00 \\
  \end{array}~
   &~
     \begin{array}{@{}c@{}@{}c@{}}
    00 \\
    01 \\
  \end{array}
    &~
      \begin{array}{@{}c@{}@{}c@{}}
    10 \\
    01 \\
  \end{array}~~
     &~
  \begin{array}{@{}c@{}@{}c@{}}
    01 \\
    10 \\
  \end{array}~
      &~
          \begin{array}{@{}c@{}@{}c@{}}
    11 \\
    01 \\
  \end{array}~
       &~
            \begin{array}{@{}c@{}@{}c@{}}
    00 \\
    10 \\
  \end{array}~
        &~
          \begin{array}{@{}c@{}@{}c@{}}
    10 \\
    10 \\
  \end{array}~
         &~
   \begin{array}{@{}c@{}@{}c@{}}
    00 \\
    11 \\
  \end{array}~
          &~
         \begin{array}{@{}c@{}@{}c@{}}
    10 \\
    11 \\
  \end{array}~~
           &~
     \begin{array}{@{}c@{}@{}c@{}}
    01 \\
    01 \\
  \end{array}~
        &~
            \begin{array}{@{}c@{}@{}c@{}}
    11 \\
    01 \\
  \end{array}~
         &~
             \begin{array}{@{}c@{}@{}c@{}}
    01 \\
    11 \\
  \end{array}~
          &~
              \begin{array}{@{}c@{}@{}c@{}}
    11 \\
    11 \\
  \end{array} ~~    \\
  \hline

    \begin{array}{@{}c@{}@{}c@{}}
    01 \\
    00 \\
  \end{array}
& ~0 & ~1 & ~1 &   ~1 &  ~0 &  0 & ~\color{red}{0} & ~\color{red}{1} & ~\color{red}{0} & ~\color{red}{1} &  ~0 &  1 &  ~1 &  ~1 &  ~1 &  0  & $A1$       \\
  \begin{array}{cc}
    11 \\
    00 \\
  \end{array}

 &   ~1  &  ~0  &  ~1  &  ~1  &  ~0  &  0  & ~\color{red}{1}  & ~\color{red}{0}  & ~\color{red}{1}  & ~\color{red}{0} &  ~1  &  0  &   ~1 &   ~1  &  ~0   & 1 & $B2$              \\
   \begin{array}{cc}
    00 \\
    00 \\
  \end{array}
  &    ~1  &  ~1  &  ~0  &  ~1  &  ~1  &  1  &  ~\color{red}{0}  &  ~\color{red}{1} &  ~1  &  ~0 &  ~\color{red}{0}  &  \color{red}{1}  &   ~0  &  ~0  &  ~0  &  1 & $E5$            \\
    \begin{array}{cc}
    10 \\
    00 \\
  \end{array}
   &    ~1  &  ~1  &  ~1  &  ~0  &  ~1  &  1  &  ~\color{red}{1}  & ~\color{red}{0} & ~0  &  ~1  &  ~\color{red}{1}  &  \color{red}{0}  &   ~0  &  ~0  &  ~1  &  0 & $F6$           \\
     \begin{array}{cc}
    00 \\
    01 \\
  \end{array}
    &  ~0  &  ~0  &  ~1  &  ~1  &  ~0  &  1  &  ~0  &  ~1  & ~\color{red}{0}  &  ~\color{red}{1}  &  ~\color{red}{0}  &  \color{red}{1}  &   ~1  &  ~1  &  ~1  &  0  & $A2$           \\
      \begin{array}{cc}
    10 \\
    01 \\
  \end{array}
     &  ~0  &  ~0  &  ~1  &  ~1  &  ~1  &  0  &  ~1  &  ~0  &  ~\color{red}{1}  &  ~\color{red}{0}  &  ~\color{red}{1}  &  \color{red}{0}  &  ~1  &  ~1  &  ~0  &  1 & $B1$            \\
     \hline
  \begin{array}{cc}
    01 \\
    10 \\
  \end{array}
      &   ~\color{red}{0}  &  ~\color{red}{1}  &  ~\color{red}{0}  &  ~\color{red}{1}  &  ~0  &  1  &  ~0  &  ~1  &  ~1  &  ~1  &  ~0  &  0  &  ~1  &  ~0  &  ~1  &  1 & $C3$          \\
          \begin{array}{cc}
    11 \\
    01 \\
  \end{array}
       &  ~\color{red}{1}  &  ~\color{red}{0}  &  ~\color{red}{1}  &  ~\color{red}{0}  &  ~1  &  0  &  ~1  &  ~0  &  ~1  &  ~1  &  ~0  &  0  &  ~0  &  ~1  &  ~1  &  1 & $D4$          \\
            \begin{array}{cc}
    00 \\
    10 \\
  \end{array}
        &   ~\color{red}{0}  &  ~\color{red}{1}  &  ~1  &  ~0  &  ~\color{red}{0}  &  \color{red}{1}  &  ~1  &  ~1  &  ~0  &  ~1  &  ~1  &  1  &  ~0  &  ~1  &  ~0  &  0  & $G7$           \\
          \begin{array}{cc}
    10 \\
    10 \\
  \end{array}
         &  ~\color{red}{1}  &  ~\color{red}{0}  &  ~0  &  ~1  &  ~\color{red}{1}  &  \color{red}{0}  &  ~1  &  ~1  &  ~1  &  ~0  &  ~1  &  1  &  ~1  &  ~0  &  ~0  &  0  & $H8$              \\
   \begin{array}{cc}
    00 \\
    11 \\
  \end{array}
          &  ~0  &  ~1  &  ~\color{red}{0}  &  ~\color{red}{1}  &  ~\color{red}{0}  &  \color{red}{1}  &  ~0  &  ~0  &  ~1  &  ~1  &  ~0  &  1  &  ~1  &  ~0  &  ~1  &  1  & $C4$             \\
         \begin{array}{cc}
    10 \\
    11 \\
  \end{array}
           & ~1  &  ~0  &  ~\color{red}{1}  &  ~\color{red}{0}  &  ~\color{red}{1}  &  \color{red}{0}  &  ~0  &  ~0  &  ~1  &  ~1  &  ~1  &  0  &  ~0  &  ~1  &  ~1  &  1 & $D3$            \\
     \hline
     \begin{array}{cc}
    01 \\
    01 \\
  \end{array}
        &  ~1  &  ~1  &  ~0  &  ~0  &  ~1  &  1  &  ~1  &  ~0  &  ~0  &  ~1  &  ~1  &  0  &   ~0  &  ~1  &  ~1  &  0 &  $F5$             \\
            \begin{array}{cc}
    11 \\
    01 \\
  \end{array}
         &  ~1  &  ~1  &  ~0  &  ~0  &  ~1  &  1  &  ~0  &  ~1  &  ~1  &  ~0  &  ~0  &  1  &  ~1  &  ~0  &  ~0  &  1 & $E6$               \\
             \begin{array}{cc}
    01 \\
    11 \\
  \end{array}
          &   ~1  &  ~0  &  ~0  &  ~1  &  ~1  &  0  &  ~1  &  ~1  &  ~0  &  ~0  &  ~1  &  1  &  ~1  &  ~0  &  ~0  &  1 & $H7$              \\
              \begin{array}{cc}
    11 \\
    11 \\
  \end{array}
  & ~0  &  ~1  &  ~1  &  ~0  &  ~0  &  1  &  ~1  &  ~1  &  ~0  &  ~0  &  ~1  &  1  &  ~0  &  ~1  &  ~1  &  0 & $G8$   \\

  \end{array}
$$
\caption*{Figure 2: The adjacency matrix, $A_{2,1}$, of the graph $\Gamma _{2,1}$}
\end{table}
The notation in the right column of Figure 2 means the following. Two rows have the same letter if and only if they correspond to non-adjacent
vertices having $8$ common neighbours; two rows have the same number if and only if they correspond to non-adjacent vertices
having $4$ common neighbours. Otherwise, every two non-adjacent vertices have $6$ common neighbours;
every two adjacent vertices have $4$ common neighbours.

\pagebreak{}
% \medskip
\noindent\textbf{5.2. The first generalisation of the smallest strictly Neumaier graph}
\medskip

In this subsection we generalise the construction from Subsection 5.1.

Take the $(e-1)$-dimensional subspace
$$
W =
\left(
  \begin{array}{ccc|cc}
  \ast &\ldots & \ast & \ast & 0 \\
  0 &\ldots &  0 & 0  & 0 \\
  \end{array}
\right),
$$
where the size of matrices is $2 \times e$.
According to Lemma \ref{PairOfGenerators}, the subspace $W$ is contained in the two generators
$$
W_1
=
\left(
  \begin{array}{ccc|cc}
   \ast &\ldots & \ast & \ast & \ast \\
    0 &\ldots & 0 & 0  & 0 \\
  \end{array}
\right)
\text{ and }
W_2
=
\left(
  \begin{array}{ccc|cc}
   \ast &\ldots & \ast & \ast & 0 \\
    0 &\ldots & 0 & 0  & \ast \\
  \end{array}
\right).
$$
Take the vector
$$
v =
\left(
  \begin{array}{ccc|cc}
    0 &\ldots & 0 & 0 & 0\\
    0 &\ldots & 0 & 1  & 0  \\
  \end{array}
\right)
$$
and consider the cosets
$$
v + W_1
=
\left(
  \begin{array}{ccc|cc}
   \ast &\ldots & \ast & \ast & \ast \\
    0 &\ldots & 0 & 1  & 0 \\
  \end{array}
\right),~~~
v + W_2
=
\left(
  \begin{array}{ccc|cc}
   \ast &\ldots & \ast & \ast & 0 \\
    0 &\ldots & 0 & 1  & \ast \\
  \end{array}
\right),
$$
whose intersection is
$$
v + W
=
\left(
  \begin{array}{ccc|cc}
   \ast &\ldots & \ast & \ast & 0 \\
    0 &\ldots & 0 & 1  & 0 \\
  \end{array}
\right).
$$

%\begin{lemma}
%The sets $W_2$ and $v+W_2$ induce $2^{e-1}$-regular $2^e$-cliques in $\Gamma(W,W_1,v)$.
%\end{lemma}

Denote by $\Gamma_{e,1}=\Gamma_{e}(W,W_1,W_2,v)$ the graph obtained from
$\Gamma_{e}=VO^+(2e,2)$ by switching edges between the cliques $W_1$, $v+W_1$ and then between the cliques $W_2$, $v+W_2$. Let $(n,k,\lambda,\mu)$ be the parameters of the affine polar graph $\Gamma_{e}=VO^+(2e,2)$ as a strongly regular graph.

\begin{theorem}\label{thm}
The graph $\Gamma_{e,1}$ is a strictly Neumaier graph with parameters $$(2^{2e},(2^{e-1}+1)(2^e-1), 2(2^{e-2}+1)(2^{e-1}-1); 2^{e-1},2^e).$$ Further, the number of common neighbours of two non-adjacent vertices in the graph takes the values $\mu - 2^{e-1},\mu$ and $\mu + 2^{e-1}$.
\end{theorem}
\proof
For any $a,b,c,d \in \mathbb{F}_2$, let
$$
  \begin{array}{c}
   \bf{ab} \\
    \bf{cd} \\
  \end{array}
$$
denote the set of matrices
$$
\left(
  \begin{array}{ccc|cc}
   \ast &\ldots & \ast & a & b \\
    0 &\ldots & 0 & c  & d \\
  \end{array}
\right).
$$
For the affine polar graph $\Gamma_{e}=VO^+(2e,2)$, consider the subgraph $\Delta$ induced by the set of all matrices
$$
\left(
  \begin{array}{ccc|cc}
   \ast &\ldots & \ast & a & b \\
    0 &\ldots & 0 & c  & d \\
  \end{array}
\right),
$$
where $a,b,c,d$ run over $\mathbb{F}_2$.
The adjacency matrix of the subgraph $\Delta$ is presented by the block-matrix in Figure 3, where
$K$ denotes the adjacency matrix of the complete graph on $2^{e-2}$ vertices;
$J$ denotes the all-ones matrix of size $2^{e-2}\times 2^{e-2}$; $Z$ denotes the all-zeroes matrix of size $2^{e-2}\times 2^{e-2}$.

\begin{table}[H]
$$
  \begin{array}{@{}c@{}|@{}c@{}@{}c@{}@{}c@{}@{}c@{}@{}c@{}@{}c@{}|@{}c@{}@{}c@{}@{}c@{}@{}c@{}@{}c@{}@{}c@{}|@{}c@{}@{}c@{}@{}c@{}@{}c@{}@{}c@{}}
  &~
  \begin{array}{@{}c@{}@{}c@{}}
    \bf{01} \\
    \bf{00} \\
  \end{array}~
&~
  \begin{array}{@{}c@{}@{}c@{}}
    \bf{11} \\
    \bf{00} \\
  \end{array}~

 &~
   \begin{array}{@{}c@{}@{}c@{}}
    \bf{00} \\
    \bf{00} \\
  \end{array}~
  &~
    \begin{array}{@{}c@{}@{}c@{}}
    \bf{10} \\
    \bf{00} \\
  \end{array}~
   &~
     \begin{array}{@{}c@{}@{}c@{}}
    \bf{00} \\
   \bf{01} \\
  \end{array}~
    &~
      \begin{array}{@{}c@{}@{}c@{}}
    \bf{10} \\
    \bf{01} \\
  \end{array}~~
     &~
  \begin{array}{@{}c@{}@{}c@{}}
    \bf{01} \\
    \bf{10} \\
  \end{array}~
      &~
          \begin{array}{@{}c@{}@{}c@{}}
    \bf{11} \\
    \bf{01} \\
  \end{array}~
       &~
            \begin{array}{@{}c@{}@{}c@{}}
    \bf{00} \\
    \bf{10} \\
  \end{array}~
        &~
          \begin{array}{@{}c@{}@{}c@{}}
    \bf{10} \\
    \bf{10} \\
  \end{array}~
         &~
   \begin{array}{@{}c@{}@{}c@{}}
    \bf{00} \\
    \bf{11} \\
  \end{array}~
          &~
         \begin{array}{@{}c@{}@{}c@{}}
    \bf{10} \\
    \bf{11} \\
  \end{array}~~
           &~
     \begin{array}{@{}c@{}@{}c@{}}
    \bf{01} \\
    \bf{01} \\
  \end{array}~
        &~
            \begin{array}{@{}c@{}@{}c@{}}
    \bf{11} \\
    \bf{01} \\
  \end{array}~
         &~
             \begin{array}{@{}c@{}@{}c@{}}
    \bf{01} \\
    \bf{11} \\
  \end{array}~
          &~
              \begin{array}{@{}c@{}@{}c@{}}
    \bf{11} \\
    \bf{11} \\
  \end{array}~     \\
  \hline

    \begin{array}{cc}
    \bf{01} \\
    \bf{00} \\
  \end{array}
& ~K & ~J & ~J &   ~J &  ~Z &  Z & ~\color{red}{J} & ~\color{red}{Z} &~\color{red}{J} & ~\color{red}{Z} &  ~Z &  J &  ~J & ~J &  ~J &  ~Z         \\
  \begin{array}{cc}
    \bf{11} \\
    \bf{00} \\
  \end{array}

 &   ~J  &  ~K  &  ~J  &  ~J &  ~Z  &  Z  & ~\color{red}{Z}  & ~\color{red}{J}  & ~\color{red}{Z}  & ~\color{red}{J} &  ~J  &  Z  &   ~J &   ~J  &  ~Z   & ~J              \\
   \begin{array}{cc}
    \bf{00} \\
    \bf{00} \\
  \end{array}
  &    ~J  &  ~J  &  ~K  &  ~J  &  ~J  &  J  &  ~\color{red}{J}  &  ~\color{red}{Z} &  ~J  &  ~Z &  ~\color{red}{J}  &  \color{red}{Z}  &   ~Z  &  ~Z  &  ~Z&  ~J           \\
    \begin{array}{cc}
    \bf{10} \\
    \bf{00} \\
  \end{array}
   &    ~J  &  ~J  &  ~J  &  ~K &  ~J  &  J  &  ~\color{red}{Z}  & ~\color{red}{J} & ~Z  &  ~J  &  ~\color{red}{Z}  &  \color{red}{J}  &   ~Z  &  ~Z  &  ~J  &  ~Z          \\
     \begin{array}{cc}
    \bf{00} \\
    \bf{01} \\
  \end{array}
    &  ~Z  &  ~Z  &  ~J  &  ~J  &  ~K  &  J  &  ~Z  &  ~J  & ~\color{red}{J}  &  ~\color{red}{Z}  &  ~\color{red}{J}  &  \color{red}{Z}  &   ~J  &  ~J  &  ~J  &  ~Z             \\
      \begin{array}{cc}
    \bf{10} \\
    \bf{01} \\
  \end{array}
     &  ~Z  &  ~Z  &  ~J  &  ~J  &  ~J  &  K  &  ~J  &  ~Z  &  ~\color{red}{Z}  &  ~\color{red}{J}  &  ~\color{red}{Z}  &  \color{red}{J}  &  ~J  &  ~J  &  ~Z  &  ~J          \\
     \hline
  \begin{array}{cc}
    \bf{01} \\
    \bf{10} \\
  \end{array}
      &   ~\color{red}{J}  &  ~\color{red}{Z}  &  ~\color{red}{J}  &  ~\color{red}{Z}  &  ~Z  &  J  &  ~K  &  ~J  &  ~J  &  ~J  &  ~Z  &  Z  &  ~J  &  ~Z &  ~J  &  ~J          \\
          \begin{array}{cc}
    \bf{11} \\
    \bf{01} \\
  \end{array}
       &  ~\color{red}{Z}  &  ~\color{red}{J}  &  ~\color{red}{Z}  &  ~\color{red}{J}  &  ~J  &  Z  &  ~J  &  ~K  &  ~J  &  ~J  &  ~Z  &  Z  &  ~Z  &  ~J  &  ~J  &  ~J          \\
            \begin{array}{cc}
    \bf{00} \\
    \bf{10} \\
  \end{array}
        &   ~\color{red}{J}  &  ~\color{red}{Z}  &  ~J  &  ~Z  &  ~\color{red}{J}  &  \color{red}{Z}  &  ~J  &  ~J  &  ~K  &  ~J  &  ~J  &  J  &  ~Z  &  ~J &  ~Z  &  ~Z             \\
          \begin{array}{cc}
    \bf{10} \\
    \bf{10} \\
  \end{array}
         &  ~\color{red}{Z}  &  ~\color{red}{J}  &  ~Z  &  ~J  &  ~\color{red}{Z}  &  \color{red}{J}  &  ~J  &  ~J  &  ~J  &  ~K  &  ~J  &  J  &  ~J  &  ~Z&  ~Z  &  ~Z               \\
   \begin{array}{cc}
    \bf{00} \\
    \bf{11} \\
  \end{array}
          &  ~Z  &  ~J  &  ~\color{red}{J}  &  ~\color{red}{Z}  &  ~\color{red}{J}  &  \color{red}{Z}  &  ~Z  &  ~Z  &  ~J  &  ~J  &  ~K  &  J  &  ~J  &  ~Z  &  ~J  &  ~J             \\
         \begin{array}{cc}
    \bf{10} \\
    \bf{11} \\
  \end{array}
           & ~J  &  ~Z  &  ~\color{red}{Z}  &  ~\color{red}{J}  &  ~\color{red}{Z}  &  \color{red}{J}  &  ~Z  &  ~Z  &  ~J  &  ~J  &  ~J  & K  &  ~Z  &  ~J &  ~J  &  ~J          \\
     \hline
     \begin{array}{cc}
    \bf{01} \\
    \bf{01} \\
  \end{array}
        &  ~J  &  ~J  &  ~Z  &  ~Z  &  ~J  &  J  &  ~J  &  ~Z  &  ~Z  &  ~J  &  ~J  &  Z  &   ~K  &  ~J  &  ~J  &  ~Z                 \\
            \begin{array}{cc}
    \bf{11} \\
    \bf{01} \\
  \end{array}
         &  ~J  &  ~J  &  ~Z  &  ~Z  &  ~J  &  J  &  ~Z  &  ~J  &  ~J  &  ~Z  &  ~Z  &  J  &  ~J  &  ~K &  ~Z  &  ~J                \\
             \begin{array}{cc}
    \bf{01} \\
    \bf{11} \\
  \end{array}
          &   ~J  &  ~Z  &  ~Z  &  ~J  &  ~J  &  Z  &  ~J  &  ~J  &  ~Z  &  ~Z  &  ~J  &  J  &  ~J  &  ~Z  &  ~K  &  ~J              \\
              \begin{array}{cc}
    \bf{11} \\
    \bf{11} \\
  \end{array}
  & ~Z  &  ~J  &  ~J  &  ~Z  &  ~Z  &  J  &  ~J  &  ~J  &  ~Z  &  ~Z  &  ~J  &  J  &  ~Z  &  ~J  &  ~J  &  ~K    \\

  \end{array}
$$
\caption*{Figure 3: The adjacency matrix, $A_{e}$, of the subgraph $\Delta$ of $\Gamma _{e}=VO^{+}(2e,2)$ }
\end{table}
Switching edges between the cliques $W_1$, $v+W_1$ and then between the cliques $W_2$, $v+W_2$
is equivalent to inverting the red entries in Figure 3. This gives the submatrix of
the adjacency matrix of $\Gamma_{e,1}=\Gamma(W,W_1,W_2,v)$ presented in Figure 4,. Note that every switched edge connects vertices from the subgraph $\Delta$.
This means that the switching preserves all edges having a vertex outside of $\Delta$.

\begin{table}[H]
$$
  \begin{array}{@{}c@{}|@{}c@{}@{}c@{}@{}c@{}@{}c@{}@{}c@{}@{}c@{}|@{}c@{}@{}c@{}@{}c@{}@{}c@{}@{}c@{}@{}c@{}|@{}c@{}@{}c@{}@{}c@{}@{}c@{}|@{}c@{}}
  & ~
  \begin{array}{@{}c@{}@{}c@{}}
    \bf{01} \\
    \bf{00} \\
  \end{array}~
&~
  \begin{array}{@{}c@{}@{}c@{}}
    \bf{11} \\
    \bf{00} \\
  \end{array}~

 &~
   \begin{array}{@{}c@{}@{}c@{}}
    \bf{00} \\
    \bf{00} \\
  \end{array}~
  &~
    \begin{array}{@{}c@{}@{}c@{}}
    \bf{10} \\
    \bf{00} \\
  \end{array}~
   &~
     \begin{array}{@{}c@{}@{}c@{}}
    \bf{00} \\
   \bf{01} \\
  \end{array}~
    &~
      \begin{array}{@{}c@{}@{}c@{}}
    \bf{10} \\
    \bf{01} \\
  \end{array}~~
     &~
  \begin{array}{@{}c@{}@{}c@{}}
    \bf{01} \\
    \bf{10} \\
  \end{array}~
      &~
          \begin{array}{@{}c@{}@{}c@{}}
    \bf{11} \\
    \bf{01} \\
  \end{array}~
       &~
            \begin{array}{@{}c@{}@{}c@{}}
    \bf{00} \\
    \bf{10} \\
  \end{array}~
        &~
          \begin{array}{@{}c@{}@{}c@{}}
    \bf{10} \\
    \bf{10} \\
  \end{array}~
         &~
   \begin{array}{@{}c@{}@{}c@{}}
    \bf{00} \\
    \bf{11} \\
  \end{array}~
          &~
         \begin{array}{@{}c@{}@{}c@{}}
    \bf{10} \\
    \bf{11} \\
  \end{array}~~
           &~
     \begin{array}{@{}c@{}@{}c@{}}
    \bf{01} \\
    \bf{01} \\
  \end{array}~
        &~
            \begin{array}{@{}c@{}@{}c@{}}
    \bf{11} \\
    \bf{01} \\
  \end{array}~
         &~
             \begin{array}{@{}c@{}@{}c@{}}
    \bf{01} \\
    \bf{11} \\
  \end{array}~
          &~
              \begin{array}{@{}c@{}@{}c@{}}
    \bf{11} \\
    \bf{11} \\
  \end{array} ~~    \\
  \hline

    \begin{array}{cc}
    \bf{01} \\
    \bf{00} \\
  \end{array}
& ~K & ~J & ~J &   ~J &  ~Z &  Z & ~\color{red}{Z} & ~\color{red}{J} &~\color{red}{Z} & ~\color{red}{J} &  ~Z &  J &  ~J & ~J &  ~J &  Z & $A1$         \\
  \begin{array}{cc}
    \bf{11} \\
    \bf{00} \\
  \end{array}

 &   ~J  &  ~K  &  ~J  &  ~J &  ~Z  &  Z  & ~\color{red}{J}  & ~\color{red}{Z}  & ~\color{red}{J}  & ~\color{red}{Z} &  ~J  &  Z  &   ~J &   ~J  &  ~Z   & J  & $B2$               \\
   \begin{array}{cc}
    \bf{00} \\
    \bf{00} \\
  \end{array}
  &    ~J  &  ~J  &  ~K  &  ~J  &  ~J  &  J  &  ~\color{red}{Z}  &  ~\color{red}{J} &  ~J  &  ~Z &  ~\color{red}{Z}  &  \color{red}{J}  &   ~Z  &  ~Z  &  ~Z&  J  & $E5$            \\
    \begin{array}{cc}
    \bf{10} \\
    \bf{00} \\
  \end{array}
   &    ~J  &  ~J  &  ~J  &  ~K &  ~J  &  J  &  ~\color{red}{J}  & ~\color{red}{Z} & ~Z  &  ~J  &  ~\color{red}{J}  &  \color{red}{Z}  &   ~Z  &  ~Z  &  ~J  &  Z  & $F6$            \\
     \begin{array}{cc}
    \bf{00} \\
    \bf{01} \\
  \end{array}
    &  ~Z  &  ~Z  &  ~J  &  ~J  &  ~K  &  J  &  ~Z  &  ~J  & ~\color{red}{Z}  &  ~\color{red}{J}  &  ~\color{red}{Z}  &  \color{red}{J}  &   ~J  &  ~J  &  ~J  &  Z   & $A2$           \\
      \begin{array}{cc}
    \bf{10} \\
    \bf{01} \\
  \end{array}
     &  ~Z  &  ~Z  &  ~J  &  ~J  &  ~J  &  K  &  ~J  &  ~Z  &  ~\color{red}{J}  &  ~\color{red}{Z}  &  ~\color{red}{J}  &  \color{red}{Z}  &  ~J  &  ~J  &  ~Z  &  J  & $B1$            \\
     \hline
  \begin{array}{cc}
    \bf{01} \\
    \bf{10} \\
  \end{array}
      &   ~\color{red}{Z}  &  ~\color{red}{J}  &  ~\color{red}{Z}  &  ~\color{red}{J}  &  ~Z  &  J  &  ~K  &  ~J  &  ~J  &  ~J  &  ~Z  &  Z  &  ~J  &  ~Z &  ~J  &  J  & $C3$          \\
          \begin{array}{cc}
    \bf{11} \\
    \bf{01} \\
  \end{array}
       &  ~\color{red}{J}  &  ~\color{red}{Z}  &  ~\color{red}{J}  &  ~\color{red}{Z}  &  ~J  &  Z  &  ~J  &  ~K  &  ~J  &  ~J  &  ~Z  &  Z  &  ~Z  &  ~J  &  ~J  &  J &  $D4$           \\
            \begin{array}{cc}
    \bf{00} \\
    \bf{10} \\
  \end{array}
        &   ~\color{red}{Z}  &  ~\color{red}{J}  &  ~J  &  ~Z  &  ~\color{red}{Z}  &  \color{red}{J}  &  ~J  &  ~J  &  ~K  &  ~J  &  ~J  &  J  &  ~Z  &  ~J &  ~Z  &  Z & $G7$            \\
          \begin{array}{cc}
    \bf{10} \\
    \bf{10} \\
  \end{array}
         &  ~\color{red}{J}  &  ~\color{red}{Z}  &  ~Z  &  ~J  &  ~\color{red}{J}  &  \color{red}{Z}  &  ~J  &  ~J  &  ~J  &  ~K  &  ~J  &  J  &  ~J  &  ~Z&  ~Z  &  Z  & $H8$               \\
   \begin{array}{cc}
    \bf{00} \\
    \bf{11} \\
  \end{array}
          &  ~Z  &  ~J  &  ~\color{red}{Z}  &  ~\color{red}{J}  &  ~\color{red}{Z}  &  \color{red}{J}  &  ~Z  &  ~Z  &  ~J  &  ~J  &  ~K  &  J  &  ~J  &  ~Z  &  ~J  &  J & $C4$               \\
         \begin{array}{cc}
    \bf{10} \\
    \bf{11} \\
  \end{array}
           & ~J  &  ~Z  &  ~\color{red}{J}  &  ~\color{red}{Z}  &  ~\color{red}{J}  &  \color{red}{Z}  &  ~Z  &  ~Z  &  ~J  &  ~J  &  ~J  & K  &  ~Z  &  ~J &  ~J  &  J & $D3$           \\
     \hline
     \begin{array}{cc}
    \bf{01} \\
    \bf{01} \\
  \end{array}
        &  ~J  &  ~J  &  ~Z  &  ~Z  &  ~J  &  J  &  ~J  &  ~Z  &  ~Z  &  ~J  &  ~J  &  Z  &   ~K  &  ~J  &  ~J  &  Z & $F5$                 \\
            \begin{array}{cc}
    \bf{11} \\
    \bf{01} \\
  \end{array}
         &  ~J  &  ~J  &  ~Z  &  ~Z  &  ~J  &  J  &  ~Z  &  ~J  &  ~J  &  ~Z  &  ~Z  &  J  &  ~J  &  ~K &  ~Z  &  J & $E6$               \\
             \begin{array}{cc}
    \bf{01} \\
    \bf{11} \\
  \end{array}
          &   ~J  &  ~Z  &  ~Z  &  ~J  &  ~J  &  Z  &  ~J  &  ~J  &  ~Z  &  ~Z  &  ~J  &  J  &  ~J  &  ~Z  &  ~K  &  J  & $H7$              \\
              \begin{array}{cc}
    \bf{11} \\
    \bf{11} \\
  \end{array}
  & ~Z  &  ~J  &  ~J  &  ~Z  &  ~Z  &  J  &  ~J  &  ~J  &  ~Z  &  ~Z  &  ~J  &  J  &  ~Z  &  ~J  &  ~J  &  K & $G8$    \\

  \end{array}
$$
\caption*{Figure 4: The adjacency matrix, $A_{e,1}$, of the subgraph $\Delta$ of $\Gamma _{e,1}=\Gamma _{e}(W,W_{1},W_{2},v)$}
\end{table}

Let $(n,k,\lambda,\mu)$ be the parameters of the affine polar graph $\Gamma_{e}=VO^+(2e,2)$ as a strongly regular graph.
We have to check that the obtained graph is a strictly Neumaier graph. Note that $W_1$ is a regular clique in $\Gamma_{e,1}=\Gamma(W,W_1,W_2,v)$.
Let us check that any pair of vertices in $\Gamma_{e,1}$ is OK,
i.e. any two adjacent vertices have $\lambda$ common neighbours.
Also, we investigate which values of $\mu$ occur in $\Gamma_{e,1}$.

Let us consider any two vertices inside of $\Delta$.
The notation in the right column of the matrix in Figure 4 means the following.
Two block-rows have the same letter if and only if any row from the one block-row
 and any row from the other block-row correspond to non-adjacent
vertices having $\mu+2^{e-1}$ common neighbours; two block-rows have the same number if and only if any row from the one block-row
 and any row from the other block-row correspond to non-adjacent
vertices having $\mu-2^{e-1}$ common neighbours. Otherwise, every two non-adjacent vertices corresponding to rows of this submatrix have $\mu$ common neighbours; every two adjacent vertices have $\lambda$ common neighbours. This means that
all pairs of vertices inside of $\Delta$ are OK.

Let us consider any two vertices outside of $\Delta$.
Their neighbours and, consequently, their common neighbours are preserved by the switching.
This means that all pairs of vertices outside of $\Delta$ are OK.

Let us consider a vertex $x$ in $\Delta$ and a vertex $y$ outside of $\Delta$.
If the neighbours of $x$ are preserved by the switching, then $x$,$y$ are OK.
Assume that the neighbours of $x$ are switched. Then the vertices $x$,$y$ are OK
since the vertex $y$ is adjacent to half the of vertices of each block of $\Delta$.
In fact, the vertex $y$ is presented by a matrix
$$
\left(
  \begin{array}{ccc|cc}
    y_1 & \ldots &y_{2e-5} & y_{2e-3} & y_{2e-1} \\
    y_2  & \ldots & y_{2e-4} & y_{2e-2} & y_{2e} \\
  \end{array}
\right),
$$
where there is at least one non-zero among $y_2, y_4, \ldots, y_{2e-4}$.
Without losing of generality, assume that $y_2 = 1$.
Let us show that $y$ is adjacent to half the of vertices in a block
$$
\left(
  \begin{array}{ccc|cc}
   \ast &\ldots & \ast & a & b \\
    0 &\ldots & 0 & c  & d \\
  \end{array}
\right).
$$
We have
$$
y +
\left(
  \begin{array}{ccc|cc}
   \ast &\ldots & \ast & a & b \\
    0 &\ldots & 0 & c  & d \\
  \end{array}
\right)
=
\left(
  \begin{array}{ccc|cc}
   \ast &\ldots & \ast & a' & b' \\
    1 &\ldots & y_{2e-4} & c'  & d' \\
  \end{array}
\right)
=
$$
$$
=
\left(
  \begin{array}{ccc|cc}
   0 &\ldots & \ast & a' & b' \\
    1 &\ldots & y_{2e-4} & c'  & d' \\
  \end{array}
\right)
\bigcup
\left(
  \begin{array}{ccc|cc}
   1 &\ldots & \ast & a' & b' \\
    1 &\ldots & y_{2e-4} & c'  & d' \\
  \end{array}
\right) = Y_0 \cup Y_1.
$$
Note that $|Y_0| = |Y_1|$, and the form $Q$ has value $0$ on one of the sets $Y_0,Y_1$ and value $1$ on the other.
We have proved that the switching preserves the number of common neighbours $x$ and $y$, completing the proof of the theorem. $\square$

\medskip
\noindent\textbf{5.3. The second construction of the smallest strictly Neumaier graph}
\medskip

Consider the graph $\Gamma_{2}=VO^+(4,2)$.
Take the generator
$$
W_1 =
\left(
  \begin{array}{cccc}
    \ast & \ast \\
    0  & 0 \\
  \end{array}
\right),
$$
the vector
$$ v =
\left(
  \begin{array}{cccc}
    0  & 0 \\
    0  & 1 \\
  \end{array}
\right)$$
and the coset
$$
v + W_1 =
\left(
  \begin{array}{cccc}
    \ast & \ast \\
    0  & 1 \\
  \end{array}
\right).
$$
Divide vertices of the $2$-regular $4$-cliques $W_1$ and $v + W_1$ into two parts as
$$
W_1 = V_0 \cup V_1,$$$$
v+ W_1 = V_2 \cup V_3,
$$
where
$$
V_0 =
\left(
  \begin{array}{ccccc}
   \ast  & 0 \\
    0   & 0 \\
  \end{array}
\right),
$$
$$
V_1 =
\left(
  \begin{array}{ccccc}
   \ast  & 1 \\
    0   & 0 \\
  \end{array}
\right),
$$
$$
V_2 =
\left(
  \begin{array}{ccccc}
   \ast & 0 \\
    0 & 1 \\
  \end{array}
\right),
$$
$$
V_3 =
\left(
  \begin{array}{ccccc}
   \ast  & 1 \\
    0  & 1 \\
  \end{array}
\right).
$$
Note that there are all possible edges between $V_0$ and $V_2$,
there are all possible edges between $V_1$ and $V_3$, there are no edges between $V_0$ and $V_3$,
and there are no edges between $V_1$ and $V_2$. Denote by $\Gamma_{2}'$ the graph obtained from $\Gamma_{2}$
by switching edges between the cliques $W_1$ and $v + W_1$. Note that each of the sets
$V_0 \cup V_3$ and $V_1 \cup V_2$ induces a $4$-clique in $\Gamma_{2}'$.

The set
$$
C :=
\left(
  \begin{array}{cccc}
    \ast & 0 \\
    1  & \ast \\
  \end{array}
\right)
$$
induces a $2$-regular $4$-clique in the graph $\Gamma_{2}'$ as well as in $\Gamma_{2}$ since the switching between
$W_1$ and $v + W_1$ did not modify the neighbourhoods of the vertices from $C$. Moreover, $C \cap (W_1 \cup v + W_1) = \emptyset$ holds,
and any vertex from $C$ is adjacent to half of the vertices of each of the sets $V_0, V_1, V_2, V_3$. This means
that the switching between the cliques $V_1 \cup V_2$, $C$ and the switching between the cliques $V_0 \cup V_3$, $C$
preserve the regularity of $\Gamma_{2}'$. Denote by $\Gamma_{2}''$ and $\Gamma_{2}'''$ the graphs obtained from $\Gamma_{2}'$ by applying
these two switchings, respectively. One can prove that the graphs  $\Gamma_{2}''$ and $\Gamma_{2}'''$ are isomorphic
to the smallest Neumaier graph. Now we show how can the adjacency matrix of the graph $\Gamma_{2}''$ be obtained from
the adjacency matrix of  $\Gamma_{2}$.

In this setting, the adjacency matrix of the affine polar graph $\Gamma_{2}=VO^+(4,2)$ can be written as in Figure 5.

\begin{table}[H]
$$
  \begin{array}{c|@{}c@{}@{}c@{}|@{}c@{}@{}c@{}|@{}c@{}@{}c@{}|@{}c@{}@{}c@{}|@{}c@{}@{}c@{}@{}c@{}@{}c@{}|@{}c@{}@{}c@{}@{}c@{}@{}c@{}}
  & ~
  \begin{array}{@{}c@{}@{}c@{}}
    00 \\
    00 \\
  \end{array} ~
& ~
  \begin{array}{@{}c@{}@{}c@{}}
    10 \\
    00 \\
  \end{array} ~~

 &~
   \begin{array}{@{}c@{}@{}c@{}}
    01 \\
    00 \\
  \end{array}~
  & ~
    \begin{array}{@{}c@{}@{}c@{}}
    11 \\
    00 \\
  \end{array} ~~
   & ~
     \begin{array}{@{}c@{}@{}c@{}}
    00 \\
    01 \\
  \end{array} ~
    & ~
      \begin{array}{@{}c@{}@{}c@{}}
    10 \\
    01 \\
  \end{array} ~~
     & ~
  \begin{array}{@{}c@{}@{}c@{}}
    01 \\
    01 \\
  \end{array} ~
      & ~
          \begin{array}{@{}c@{}@{}c@{}}
    11 \\
    01 \\
  \end{array} ~~
       & ~
            \begin{array}{@{}c@{}@{}c@{}}
    00 \\
    10 \\
  \end{array} ~
        & ~
          \begin{array}{@{}c@{}@{}c@{}}
    10 \\
    10 \\
  \end{array} ~
         & ~
   \begin{array}{@{}c@{}@{}c@{}}
    00 \\
    11 \\
  \end{array} ~
          & ~
         \begin{array}{@{}c@{}@{}c@{}}
    10 \\
    11 \\
  \end{array} ~~
           & ~
     \begin{array}{@{}c@{}@{}c@{}}
    01 \\
    10 \\
  \end{array} ~
        & ~
            \begin{array}{@{}c@{}@{}c@{}}
    11 \\
    10 \\
  \end{array} ~
         & ~
             \begin{array}{@{}c@{}@{}c@{}}
    01 \\
    11 \\
  \end{array} ~
          & ~
              \begin{array}{@{}c@{}@{}c@{}}
    11 \\
    11 \\
  \end{array}  ~   \\
  \hline

    \begin{array}{@{}c@{}@{}c@{}}
    00 \\
    00 \\
  \end{array}
& ~0 & 1 & ~1 & 1 & ~\color{red}{1} & \color{red}{1} & ~\color{red}{0} & \color{red}{0} & ~1 & ~0 & ~1 & 0 & ~1 & ~0 & ~0 & ~1        \\
  \begin{array}{@{}c@{}@{}c@{}}
    10 \\
    00 \\
  \end{array}
 & ~1 & 0 & ~1 & 1 & ~\color{red}{1} & \color{red}{1} & ~\color{red}{0} & \color{red}{0} & ~0 & ~1 & ~0 & 1 & ~0 & ~1 & ~1 & ~0  \\
 \hline
   \begin{array}{@{}c@{}@{}c@{}}
    01 \\
    00 \\
  \end{array}
  & ~1 & 1 & ~0 & 1 & ~\color{red}{0} & \color{red}{0} & ~\color{red}{1} & \color{red}{1} & ~\color{red}{1} & ~\color{red}{0} & ~\color{red}{0} & \color{red}{1} & ~1 & ~0 & ~1 & ~0   \\
    \begin{array}{@{}c@{}@{}c@{}}
    11 \\
    00 \\
  \end{array}
   & ~1 & 1 & ~1 & 0 & ~\color{red}{0} & \color{red}{0} & ~\color{red}{1} & \color{red}{1} & ~\color{red}{0} & ~\color{red}{1} & ~\color{red}{1} & \color{red}{0} & ~0 & ~1 & ~0 & ~1   \\
   \hline
     \begin{array}{@{}c@{}@{}c@{}}
    00 \\
    01 \\
  \end{array}
    & ~\color{red}{1} & \color{red}{1} & ~\color{red}{0} & \color{red}{0} & ~0 & 1 & ~1 & 1 & ~\color{red}{1} & ~\color{red}{0} & ~\color{red}{1} & \color{red}{0} & ~0 & ~1 & ~1 & ~0  \\
      \begin{array}{@{}c@{}@{}c@{}}
    10 \\
    01 \\
  \end{array}
    & ~\color{red}{1} & \color{red}{1} & ~\color{red}{0} & \color{red}{0} & ~1 & 0 & ~1 & 1 & ~\color{red}{0} & ~\color{red}{1} & ~\color{red}{0} & \color{red}{1} & ~1 & ~0 & ~0 & ~1  \\
     \hline
  \begin{array}{@{}c@{}@{}c@{}}
    01 \\
    01 \\
  \end{array}
      & ~\color{red}{0} & \color{red}{0} & ~\color{red}{1} & \color{red}{1} & ~1 & 1 & ~0 & 1 & ~0 & ~1 & ~1 & 0 & ~1 & ~0 & ~1 & ~0      \\
          \begin{array}{@{}c@{}@{}c@{}}
    11 \\
    01 \\
  \end{array}
       & ~\color{red}{0} & \color{red}{0} & ~\color{red}{1} & \color{red}{1} & ~1 & 1 & ~1 & 0 & ~1 & ~0 & ~0 & 1 & ~0 & ~1 & ~0 & ~1     \\
       \hline
            \begin{array}{@{}c@{}@{}c@{}}
    00 \\
    10 \\
  \end{array}
       & ~1 & 0 & ~\color{red}1 & \color{red}0 & ~\color{red}1 & \color{red}0 & ~0 & 1 & ~0 & ~1 & ~1 & 1 & ~1 & ~1 & ~0 & ~0     \\
          \begin{array}{@{}c@{}@{}c@{}}
    10 \\
    10 \\
  \end{array}
       & ~0 & 1 & ~\color{red}0 & \color{red}1 & ~\color{red}0 & \color{red}1 & ~1 & 0 & ~1 & ~0 & ~1 & 1 & ~1 & ~1 & ~0 & ~0      \\
   \begin{array}{@{}c@{}@{}c@{}}
    00 \\
    11 \\
  \end{array}
      & ~1 & 0 & ~\color{red}0 & \color{red}1 & ~\color{red}1 & \color{red}0 & ~1 & 0 & ~1 & ~1 & ~0 & 1 & ~0 & ~0 & ~1 & ~1    \\
         \begin{array}{@{}c@{}@{}c@{}}
    10 \\
    11 \\
  \end{array}
     & ~0 & 1 & ~\color{red}1 & \color{red}0 & ~\color{red}0 & \color{red}1 & ~0 & 1 & ~1 & ~1 & ~1 & 0 & ~0 & ~0 & ~1 & ~1     \\
     \hline
     \begin{array}{@{}c@{}@{}c@{}}
    01 \\
    10 \\
  \end{array}
    & ~1 & 0 & ~1 & 0 & ~0 & 1 & ~1 & 0 & ~1 & ~1 & ~0 & 0 & ~0 & ~1 & ~1 & ~1     \\
            \begin{array}{@{}c@{}@{}c@{}}
    11 \\
    10 \\
  \end{array}
     & ~0 & 1 & ~0 & 1 & ~1 & 0 & ~0 & 1 & ~1 & ~1 & ~0 & 0 & ~1 & ~0 & ~1 & ~1    \\
             \begin{array}{@{}c@{}@{}c@{}}
    01 \\
    11 \\
  \end{array}
     & ~0 & 1 & ~1 & 0 & ~1 & 0 & ~1 & 0 & ~0 & ~0 & ~1 & 1 & ~1 & ~1 & ~0 & ~1   \\
              \begin{array}{@{}c@{}@{}c@{}}
    11 \\
    11 \\
  \end{array}
  & ~1 & 0 & ~0 & 1 & ~0 & 1 & ~0 & 1 & ~0 & ~0 & ~1 & 1 & ~1 & ~1 & ~1 & ~0  \\
  \end{array}
$$
\caption*{Figure 5: The adjacency matrix, $B_{2}$, of  $\Gamma_{2}$}
\end{table}
Switching edges between the cliques $W_1$, $v+W_1$ and then between the cliques $V_1 \cup V_2$, $C$
is equivalent to inverting the red entries in Figure 5. This gives the strictly Neumaier graph $\Gamma_{2,2}$ on 16 vertices,
whose adjacency matrix is presented in Figure 6.

\begin{table}[H]
$$
  \begin{array}{c|@{}c@{}@{}c@{}|@{}c@{}@{}c@{}|@{}c@{}@{}c@{}|@{}c@{}@{}c@{}|@{}c@{}@{}c@{}@{}c@{}@{}c@{}|@{}c@{}@{}c@{}@{}c@{}@{}c@{}|@{}c@{}}
  & ~
  \begin{array}{@{}c@{}@{}c@{}}
    00 \\
    00 \\
  \end{array} ~
& ~
  \begin{array}{@{}c@{}@{}c@{}}
    10 \\
    00 \\
  \end{array} ~~

 &~
   \begin{array}{@{}c@{}@{}c@{}}
    01 \\
    00 \\
  \end{array}~
  & ~
    \begin{array}{@{}c@{}@{}c@{}}
    11 \\
    00 \\
  \end{array} ~~
   & ~
     \begin{array}{@{}c@{}@{}c@{}}
    00 \\
    01 \\
  \end{array} ~
    & ~
      \begin{array}{@{}c@{}@{}c@{}}
    10 \\
    01 \\
  \end{array} ~~
     & ~
  \begin{array}{@{}c@{}@{}c@{}}
    01 \\
    01 \\
  \end{array} ~
      & ~
          \begin{array}{@{}c@{}@{}c@{}}
    11 \\
    01 \\
  \end{array} ~~
       & ~
            \begin{array}{@{}c@{}@{}c@{}}
    00 \\
    10 \\
  \end{array} ~
        & ~
          \begin{array}{@{}c@{}@{}c@{}}
    10 \\
    10 \\
  \end{array} ~
         & ~
   \begin{array}{@{}c@{}@{}c@{}}
    00 \\
    11 \\
  \end{array} ~
          & ~
         \begin{array}{@{}c@{}@{}c@{}}
    10 \\
    11 \\
  \end{array} ~~
           & ~
     \begin{array}{@{}c@{}@{}c@{}}
    01 \\
    10 \\
  \end{array} ~
        & ~
            \begin{array}{@{}c@{}@{}c@{}}
    11 \\
    10 \\
  \end{array} ~
         & ~
             \begin{array}{@{}c@{}@{}c@{}}
    01 \\
    11 \\
  \end{array} ~
          & ~
              \begin{array}{@{}c@{}@{}c@{}}
    11 \\
    11 \\
  \end{array}  ~~   \\
  \hline

    \begin{array}{@{}c@{}@{}c@{}}
    00 \\
    00 \\
  \end{array}
& ~0 & 1 & ~1 & 1 & ~\color{red}{0} & \color{red}{0} & ~\color{red}{1} & \color{red}{1} & ~1 & ~0 & ~1 & 0 & ~1 & ~0 & ~0 & 1 & $A1$  \\
  \begin{array}{@{}c@{}@{}c@{}}
    10 \\
    00 \\
  \end{array}
 & ~1 & 0 & ~1 & 1 & ~\color{red}{0} & \color{red}{0} & ~\color{red}{1} & \color{red}{1} & ~0 & ~1 & ~0 & 1 & ~0 & ~1 & ~1 & 0 & $B2$  \\
 \hline
   \begin{array}{@{}c@{}@{}c@{}}
    01 \\
    00 \\
  \end{array}
  & ~1 & 1 & ~0 & 1 & ~\color{red}{1} & \color{red}{1} & ~\color{red}{0} & \color{red}{0} & ~\color{red}{0} & ~\color{red}{1} & ~\color{red}{1} & \color{red}{0} & ~1 & ~0 & ~1 & 0 & $C3$  \\
    \begin{array}{@{}c@{}@{}c@{}}
    11 \\
    00 \\
  \end{array}
   & ~1 & 1 & ~1 & 0 & ~\color{red}{1} & \color{red}{1} & ~\color{red}{0} & \color{red}{0} & ~\color{red}{1} & ~\color{red}{0} & ~\color{red}{0} & \color{red}{1} & ~0 & ~1 & ~0 & 1 & $D4$  \\
   \hline
     \begin{array}{@{}c@{}@{}c@{}}
    00 \\
    01 \\
  \end{array}
    & ~\color{red}{0} & \color{red}{0} & ~\color{red}{1} & \color{red}{1} & ~0 & 1 & ~1 & 1 & ~\color{red}{0} & ~\color{red}{1} & ~\color{red}{0} & \color{red}{1} & ~0 & ~1 & ~1 & 0 & $B1$ \\
      \begin{array}{@{}c@{}@{}c@{}}
    10 \\
    01 \\
  \end{array}
    & ~\color{red}{0} & \color{red}{0} & ~\color{red}{1} & \color{red}{1} & ~1 & 0 & ~1 & 1 & ~\color{red}{1} & ~\color{red}{0} & ~\color{red}{1} & \color{red}{0} & ~1 & ~0 & ~0 & 1 & $A2$  \\
     \hline
  \begin{array}{@{}c@{}@{}c@{}}
    01 \\
    01 \\
  \end{array}
      & ~\color{red}{1} & \color{red}{1} & ~\color{red}{0} & \color{red}{0} & ~1 & 1 & ~0 & 1 & ~0 & ~1 & ~1 & 0 & ~1 & ~0 & ~1 & 0 & $C4$ \\
          \begin{array}{@{}c@{}@{}c@{}}
    11 \\
    01 \\
  \end{array}
       & ~\color{red}{1} & \color{red}{1} & ~\color{red}{0} & \color{red}{0} & ~1 & 1 & ~1 & 0 & ~1 & ~0 & ~0 & 1 & ~0 & ~1 & ~0 & 1 & $D3$    \\
       \hline
            \begin{array}{@{}c@{}@{}c@{}}
    00 \\
    10 \\
  \end{array}
       & ~1 & 0 & ~\color{red}{0} & \color{red}{1} & ~\color{red}{0} & \color{red}{1} & ~0 & 1 & ~0 & ~1 & ~1 & 1 & ~1 & ~1 & ~0 & 0 & $E5$     \\
          \begin{array}{@{}c@{}@{}c@{}}
    10 \\
    10 \\
  \end{array}
       & ~0 & 1 & ~\color{red}{1} & \color{red}{0} & ~\color{red}{1} & \color{red}{0} & ~1 & 0 & ~1 & ~0 & ~1 & 1 & ~1 & ~1 & ~0 & 0 & $F6$   \\
   \begin{array}{@{}c@{}@{}c@{}}
    00 \\
    11 \\
  \end{array}
      & ~1 & 0 & ~\color{red}{1} & \color{red}{0} & ~\color{red}{0} & \color{red}{1} & ~1 & 0 & ~1 & ~1 & ~0 & 1 & ~0 & ~0 & ~1 & 1 & $G7$   \\
         \begin{array}{@{}c@{}@{}c@{}}
    10 \\
    11 \\
  \end{array}
     & ~0 & 1 & ~\color{red}{0} & \color{red}{1} & ~\color{red}{1} & \color{red}{0} & ~0 & 1 & ~1 & ~1 & ~1 & 0 & ~0 & ~0 & ~1 & 1& $H8$     \\
     \hline
     \begin{array}{@{}c@{}@{}c@{}}
    01 \\
    10 \\
  \end{array}
    & ~1 & 0 & ~1 & 0 & ~0 & 1 & ~1 & 0 & ~1 & ~1 & ~0 & 0 & ~0 & ~1 & ~1 & 1 & $G8$     \\
            \begin{array}{@{}c@{}@{}c@{}}
    11 \\
    10 \\
  \end{array}
     & ~0 & 1 & ~0 & 1 & ~1 & 0 & ~0 & 1 & ~1 & ~1 & ~0 & 0 & ~1 & ~0 & ~1 & 1 & $H7$   \\
             \begin{array}{@{}c@{}@{}c@{}}
    01 \\
    11 \\
  \end{array}
     & ~0 & 1 & ~1 & 0 & ~1 & 0 & ~1 & 0 & ~0 & ~0 & ~1 & 1 & ~1 & ~1 & ~0 & 1 & $F5$   \\
              \begin{array}{@{}c@{}@{}c@{}}
    11 \\
    11 \\
  \end{array}
  & ~1 & 0 & ~0 & 1 & ~0 & 1 & ~0 & 1 & ~0 & ~0 & ~1 & 1 & ~1 & ~1 & ~1 & 0 & $E6$  \\
  \end{array}
$$
\caption*{Figure 6: The adjacency matrix, $B_{2,2}$, of $\Gamma _{2,2}$}
\end{table}
The notation in the right column of Figure 6 means the following. Two rows have the same letter if and only if they correspond to non-adjacent
vertices having $8$ common neighbours; two rows have the same number if and only if they correspond to non-adjacent vertices
having $4$ common neighbours. Otherwise, every two non-adjacent vertices have $6$ common neighbours;
every two adjacent vertices have $4$ common neighbours.

\pagebreak{}
% \medskip
\noindent\textbf{5.4. The second generalisation of the smallest strictly Neumaier graph}
\medskip

In this subsection we generalise the construction from Subsection 5.3 and present one more family of strictly Neumaier graphs.

For any $e \ge 2$, consider the affine polar graph $\Gamma_{e}=VO^+(2e, 2)$ and
take the regular clique given by the generator
$$
W_1 =
\left(
  \begin{array}{ccc|cc}
  \ast &\ldots & \ast & \ast & \ast \\
  0 &\ldots &  0 & 0  & 0 \\
  \end{array}
\right).
$$
For the vector
$$
v =
\left(
  \begin{array}{ccc|cc}
    0 &\ldots & 0 & 0 & 0\\
    0 &\ldots & 0 & 0 & 1  \\
  \end{array}
\right),
$$
take the regular clique
$$
v + W_1
=
\left(
  \begin{array}{ccc|cc}
   \ast &\ldots & \ast & \ast & \ast \\
    0 &\ldots & 0 & 0  & 1 \\
  \end{array}
\right),
$$
which lies in the spread given by $W_1$. Divide $W_1$ and $v + W_1$ into two parts as
$$
W_1 = V_0 \cup V_1,$$$$
v+ W_1 = V_2 \cup V_3,
$$
where
$$
V_0 =
\left(
  \begin{array}{ccc|cc}
   \ast &\ldots & \ast & \ast & 0 \\
    0 &\ldots & 0 & 0  & 0 \\
  \end{array}
\right),
$$
$$
V_1 =
\left(
  \begin{array}{ccc|cc}
   \ast &\ldots & \ast & \ast & 1 \\
    0 &\ldots & 0 & 0  & 0 \\
  \end{array}
\right),
$$
$$
V_2 =
\left(
  \begin{array}{ccc|cc}
   \ast &\ldots & \ast & \ast & 0 \\
    0 &\ldots & 0 & 0  & 1 \\
  \end{array}
\right),
$$
$$
V_3 =
\left(
  \begin{array}{ccc|cc}
   \ast &\ldots & \ast & \ast & 1 \\
    0 &\ldots & 0 & 0  & 1 \\
  \end{array}
\right).
$$

Note that there are all possible edges between $V_0$ and $V_2$,
there are all possible edges between $V_1$ and $V_3$, there are no edges between $V_0$ and $V_3$,
and there are no edges between $V_1$ and $V_2$. Denote by $\Gamma_{e}'$ the graph obtained from $\Gamma_{e}$
by switching edges between the cliques $W_1$ and $v + W_1$. Note that each of the sets
$V_0 \cup V_3$ and $V_1 \cup V_2$ induces a $2^e$-clique in $\Gamma_{e}'$.

The set
$$
C :=
\left(
  \begin{array}{ccc|cc}
   \ast & \ldots & \ast & \ast & 0 \\
   0 & \ldots & 0 & 1  & \ast \\
  \end{array}
\right)
$$
induces a $2^{e-1}$-regular $2^e$-clique in the graph $\Gamma_{e}'$ as well as in $\Gamma_{e}$ since the switching between
$W_1$ and $v + W_1$ did not modify the neighbourhoods of the vertices from $C$. Moreover, $C \cap (W_1 \cup v + W_1) = \emptyset$ holds,
and any vertex from $C$ is adjacent to half of the vertices of each of the sets $V_0, V_1, V_2, V_3$. This means
that the switching between the cliques $V_1 \cup V_2$, $C$ and the switching between the cliques $V_0 \cup V_3$, $C$
preserve the regularity of $\Gamma_{e}'$.
Denote by $\Gamma_{e,2}$ the graph obtained from $\Gamma_{e}'$
by switching edges between the cliques $W_1 \cup W_2$ and $C$. Let $(n,k,\lambda,\mu)$ be the parameters of the affine polar graph $\Gamma_{e}=VO^+(2e,2)$ as a strongly regular graph.

\begin{theorem}
The graph $\Gamma_{e,2}$ is a strictly Neumaier graph with parameters
$$(2^{2e},(2^{e-1}+1)(2^e-1), 2(2^{e-2}+1)(2^{e-1}-1); 2^{e-1},2^e).$$
Further, the number of common neighbours of two non-adjacent vertices in the graph takes the values $\mu - 2^{e-1},\mu$ and $\mu + 2^{e-1}$.
\end{theorem}
\proof
For any $a,b,c,d \in \mathbb{F}_2$, let
$$
  \begin{array}{c}
   \bf{ab} \\
    \bf{cd} \\
  \end{array}
$$
denote the set of matrices
$$
\left(
  \begin{array}{ccc|cc}
   \ast &\ldots & \ast & a & b \\
    0 &\ldots & 0 & c  & d \\
  \end{array}
\right).
$$
For the affine polar graph $\Gamma_{e}=VO^+(2e,2)$, consider the subgraph $\Delta$ induced by the set of all matrices
$$
\left(
  \begin{array}{ccc|cc}
   \ast &\ldots & \ast & a & b \\
    0 &\ldots & 0 & c  & d \\
  \end{array}
\right),
$$
where $a,b,c,d$ run over $\mathbb{F}_2$.
The adjacency matrix of the subgraph $\Delta$ is presented by the block-matrix in Figure 7, where
$K$ denotes the adjacency matrix of the complete graph on $2^{e-2}$ vertices;
$J$ denotes the all-ones matrix of size $2^{e-2}\times 2^{e-2}$; $Z$ denotes the all-zeroes matrix of size $2^{e-2}\times 2^{e-2}$.

\begin{table}[H]
$$
  \begin{array}{c|@{}c@{}@{}c@{}|@{}c@{}@{}c@{}|@{}c@{}@{}c@{}|@{}c@{}@{}c@{}|@{}c@{}@{}c@{}@{}c@{}@{}c@{}|@{}c@{}@{}c@{}@{}c@{}@{}c@{}}
  & ~
  \begin{array}{@{}c@{}@{}c@{}}
    \bf{00} \\
    \bf{00} \\
  \end{array} ~
& ~
  \begin{array}{@{}c@{}@{}c@{}}
    \bf{10} \\
    \bf{00} \\
  \end{array} ~~

 &~
   \begin{array}{@{}c@{}@{}c@{}}
    \bf{01} \\
    \bf{00} \\
  \end{array}~
  & ~
    \begin{array}{@{}c@{}@{}c@{}}
    \bf{11} \\
    \bf{00} \\
  \end{array} ~~
   & ~
     \begin{array}{@{}c@{}@{}c@{}}
    \bf{00} \\
    \bf{01} \\
  \end{array} ~
    & ~
      \begin{array}{@{}c@{}@{}c@{}}
    \bf{10} \\
    \bf{01} \\
  \end{array} ~~
     & ~
  \begin{array}{@{}c@{}@{}c@{}}
    \bf{01} \\
    \bf{01} \\
  \end{array} ~
      & ~
          \begin{array}{@{}c@{}@{}c@{}}
    \bf{11} \\
    \bf{01} \\
  \end{array} ~~
       & ~
            \begin{array}{@{}c@{}@{}c@{}}
    \bf{00} \\
    \bf{10} \\
  \end{array} ~
        & ~
          \begin{array}{@{}c@{}@{}c@{}}
    \bf{10} \\
    \bf{10} \\
  \end{array} ~
         & ~
   \begin{array}{@{}c@{}@{}c@{}}
    \bf{00} \\
    \bf{11} \\
  \end{array} ~
          & ~
         \begin{array}{@{}c@{}@{}c@{}}
    \bf{10} \\
    \bf{11} \\
  \end{array} ~~
           & ~
     \begin{array}{@{}c@{}@{}c@{}}
    \bf{01} \\
    \bf{10} \\
  \end{array} ~
        & ~
            \begin{array}{@{}c@{}@{}c@{}}
    \bf{11} \\
    \bf{10} \\
  \end{array} ~
         & ~
             \begin{array}{@{}c@{}@{}c@{}}
    \bf{01} \\
    \bf{11} \\
  \end{array} ~
          & ~
              \begin{array}{@{}c@{}@{}c@{}}
    \bf{11} \\
    \bf{11} \\
  \end{array}  ~   \\
  \hline

    \begin{array}{@{}c@{}@{}c@{}}
    \bf{00} \\
    \bf{00} \\
  \end{array}
& ~K & J & ~J & J & ~\color{red}{J} & \color{red}{J} & ~\color{red}{Z} & \color{red}{Z} & ~J & ~Z & ~J & Z & ~J & ~Z & ~Z & ~J        \\
  \begin{array}{@{}c@{}@{}c@{}}
    \bf{10} \\
    \bf{00} \\
  \end{array}
 & ~J & K & ~J & J & ~\color{red}{J} & \color{red}{J} & ~\color{red}{Z} & \color{red}{Z} & ~Z & ~J & ~Z & J & ~Z & ~J & ~J & ~Z  \\
 \hline
   \begin{array}{@{}c@{}@{}c@{}}
    \bf{01} \\
    \bf{00} \\
  \end{array}
  & ~J & J & ~K & J & ~\color{red}{Z} & \color{red}{Z} & ~\color{red}{J} & \color{red}{J} & ~\color{red}{J} & ~\color{red}{Z} & ~\color{red}{Z} & \color{red}{J} & ~J & ~Z & ~J & ~Z   \\
    \begin{array}{@{}c@{}@{}c@{}}
    \bf{11} \\
    \bf{00} \\
  \end{array}
   & ~J & J & ~J & K & ~\color{red}{Z} & \color{red}{Z} & ~\color{red}{J} & \color{red}{J} & ~\color{red}{Z} & ~\color{red}{J} & ~\color{red}{J} & \color{red}{Z} & ~Z & ~J & ~Z & ~J   \\
   \hline
     \begin{array}{@{}c@{}@{}c@{}}
    \bf{00} \\
    \bf{01} \\
  \end{array}
    & ~\color{red}{J} & \color{red}{J} & ~\color{red}{Z} & \color{red}{Z} & ~K & J & ~J & J & ~\color{red}{J} & ~\color{red}{Z} & ~\color{red}{J} & \color{red}{Z} & ~Z & ~J & ~J & ~Z  \\
      \begin{array}{@{}c@{}@{}c@{}}
    \bf{10} \\
    \bf{01} \\
  \end{array}
    & ~\color{red}{J} & \color{red}{J} & ~\color{red}{Z} & \color{red}{Z} & ~J & K & ~J & J & ~\color{red}{Z} & ~\color{red}{J} & ~\color{red}{Z} & \color{red}{J} & ~J & ~Z & ~Z & ~J  \\
     \hline
  \begin{array}{@{}c@{}@{}c@{}}
    \bf{01} \\
    \bf{01} \\
  \end{array}
      & ~\color{red}{Z} & \color{red}{Z} & ~\color{red}{J} & \color{red}{J} & ~J & J & ~K & J & ~Z & ~J & ~J & Z & ~J & ~Z & ~J & ~Z      \\
          \begin{array}{@{}c@{}@{}c@{}}
    \bf{11} \\
    \bf{01} \\
  \end{array}
       & ~\color{red}{Z} & \color{red}{Z} & ~\color{red}{J} & \color{red}{J} & ~J & J & ~J & K & ~J & ~Z & ~Z & J & ~Z & ~J & ~Z & ~J     \\
       \hline
            \begin{array}{@{}c@{}@{}c@{}}
    \bf{00} \\
    \bf{10} \\
  \end{array}
       & ~J & Z & ~\color{red}J & \color{red}Z & ~\color{red}J & \color{red}Z & ~Z & J & ~K & ~J & ~J & J & ~J & ~J & ~Z & ~Z     \\
          \begin{array}{@{}c@{}@{}c@{}}
    \bf{10} \\
    \bf{10} \\
  \end{array}
       & ~Z & J & ~\color{red}Z & \color{red}J & ~\color{red}Z & \color{red}J & ~J & Z & ~J & ~K & ~J & J & ~J & ~J & ~Z & ~Z      \\
   \begin{array}{@{}c@{}@{}c@{}}
    \bf{00} \\
    \bf{11} \\
  \end{array}
      & ~J & Z & ~\color{red}Z & \color{red}J & ~\color{red}J & \color{red}Z & ~J & Z & ~J & ~J & ~K & J & ~Z & ~Z & ~J & ~J    \\
         \begin{array}{@{}c@{}@{}c@{}}
    \bf{10} \\
    \bf{11} \\
  \end{array}
     & ~Z & J & ~\color{red}J & \color{red}Z & ~\color{red}Z & \color{red}J & ~Z & J & ~J & ~J & ~J & K & ~Z & ~Z & ~J & ~J     \\
     \hline
     \begin{array}{@{}c@{}@{}c@{}}
    \bf{01} \\
    \bf{10} \\
  \end{array}
    & ~J & Z & ~J & Z & ~Z & J & ~J & Z & ~J & ~J & ~Z & Z & ~K & ~J & ~J & ~J     \\
            \begin{array}{@{}c@{}@{}c@{}}
    \bf{11} \\
    \bf{10} \\
  \end{array}
     & ~Z & J & ~Z & J & ~J & Z & ~Z & J & ~J & ~J & ~Z & Z & ~J & ~K & ~J & ~J    \\
             \begin{array}{@{}c@{}@{}c@{}}
    \bf{01} \\
    \bf{11} \\
  \end{array}
     & ~Z & J & ~J & Z & ~J & Z & ~J & Z & ~Z & ~Z & ~J & J & ~J & ~J & ~K & ~J   \\
              \begin{array}{@{}c@{}@{}c@{}}
    \bf{11} \\
    \bf{11} \\
  \end{array}
  & ~J & Z & ~Z & J & ~Z & J & ~Z & J & ~Z & ~Z & ~J & J & ~J & ~J & ~J & ~K  \\
  \end{array}
$$
\caption*{Figure 7: The adjacency matrix, $B_{e}$, of the subgraph $\Delta$ of $\Gamma _{e}$}
\end{table}
Switching edges between the cliques $W_1$, $v+W_1$ and then between the cliques $V_1 \cup V_2$, $C$
is equivalent to inverting the red entries in Figure 7. This gives the submatrix of
the adjacency matrix of $\Gamma_{e,2}$ presented in Figure 8. Note that every switched edge connects vertices from the subgraph $\Delta$.
This means that the switching preserves all edges having a vertex outside of $\Delta$.

\begin{table}[H]
$$
  \begin{array}{c|@{}c@{}@{}c@{}|@{}c@{}@{}c@{}|@{}c@{}@{}c@{}|@{}c@{}@{}c@{}|@{}c@{}@{}c@{}@{}c@{}@{}c@{}|@{}c@{}@{}c@{}@{}c@{}@{}c@{}|@{}c@{}}
  & ~
  \begin{array}{@{}c@{}@{}c@{}}
    \bf{00} \\
    \bf{00} \\
  \end{array} ~
& ~
  \begin{array}{@{}c@{}@{}c@{}}
    \bf{10} \\
    \bf{00} \\
  \end{array} ~~

 &~
   \begin{array}{@{}c@{}@{}c@{}}
    \bf{01} \\
    \bf{00} \\
  \end{array}~
  & ~
    \begin{array}{@{}c@{}@{}c@{}}
    \bf{11} \\
    \bf{00} \\
  \end{array} ~~
   & ~
     \begin{array}{@{}c@{}@{}c@{}}
    \bf{00} \\
    \bf{01} \\
  \end{array} ~
    & ~
      \begin{array}{@{}c@{}@{}c@{}}
    \bf{10} \\
    \bf{01} \\
  \end{array} ~~
     & ~
  \begin{array}{@{}c@{}@{}c@{}}
    \bf{01} \\
    \bf{01} \\
  \end{array} ~
      & ~
          \begin{array}{@{}c@{}@{}c@{}}
    \bf{11} \\
    \bf{01} \\
  \end{array} ~~
       & ~
            \begin{array}{@{}c@{}@{}c@{}}
    \bf{00} \\
    \bf{10} \\
  \end{array} ~
        & ~
          \begin{array}{@{}c@{}@{}c@{}}
    \bf{10} \\
    \bf{10} \\
  \end{array} ~
         & ~
   \begin{array}{@{}c@{}@{}c@{}}
    \bf{00} \\
    \bf{11} \\
  \end{array} ~
          & ~
         \begin{array}{@{}c@{}@{}c@{}}
    \bf{10} \\
    \bf{11} \\
  \end{array} ~~
           & ~
     \begin{array}{@{}c@{}@{}c@{}}
    \bf{01} \\
    \bf{10} \\
  \end{array} ~
        & ~
            \begin{array}{@{}c@{}@{}c@{}}
    \bf{11} \\
    \bf{10} \\
  \end{array} ~
         & ~
             \begin{array}{@{}c@{}@{}c@{}}
    \bf{01} \\
    \bf{11} \\
  \end{array} ~
          & ~
              \begin{array}{@{}c@{}@{}c@{}}
    \bf{11} \\
    \bf{11} \\
  \end{array}  ~~   \\
  \hline

    \begin{array}{@{}c@{}@{}c@{}}
    \bf{00} \\
    \bf{00} \\
  \end{array}
& ~K & J & ~J & J & ~\color{red}{Z} & \color{red}{Z} & ~\color{red}{J} & \color{red}{J} & ~J & ~Z & ~J & Z & ~J & ~Z & ~Z & J & $A1$  \\
  \begin{array}{@{}c@{}@{}c@{}}
    \bf{10} \\
    \bf{00} \\
  \end{array}
 & ~J & K & ~J & J & ~\color{red}{Z} & \color{red}{Z} & ~\color{red}{J} & \color{red}{J} & ~Z & ~J & ~Z & J & ~Z & ~J & ~J & Z & $B2$  \\
 \hline
   \begin{array}{@{}c@{}@{}c@{}}
    \bf{01} \\
    \bf{00} \\
  \end{array}
  & ~J & J & ~K & J & ~\color{red}{J} & \color{red}{J} & ~\color{red}{Z} & \color{red}{Z} & ~\color{red}{Z} & ~\color{red}{J} & ~\color{red}{J} & \color{red}{Z} & ~J & ~Z & ~J & Z & $C3$  \\
    \begin{array}{@{}c@{}@{}c@{}}
    \bf{11} \\
    \bf{00} \\
  \end{array}
   & ~J & J & ~J & K & ~\color{red}{J} & \color{red}{J} & ~\color{red}{Z} & \color{red}{Z} & ~\color{red}{J} & ~\color{red}{Z} & ~\color{red}{Z} & \color{red}{J} & ~Z & ~J & ~Z & J & $D4$  \\
   \hline
     \begin{array}{@{}c@{}@{}c@{}}
    \bf{00} \\
    \bf{01} \\
  \end{array}
    & ~\color{red}{Z} & \color{red}{Z} & ~\color{red}{J} & \color{red}{J} & ~K & J & ~J & J & ~\color{red}{Z} & ~\color{red}{J} & ~\color{red}{Z} & \color{red}{J} & ~Z & ~J & ~J & Z & $B1$ \\
      \begin{array}{@{}c@{}@{}c@{}}
    \bf{10} \\
    \bf{01} \\
  \end{array}
    & ~\color{red}{Z} & \color{red}{Z} & ~\color{red}{J} & \color{red}{J} & ~J & K & ~J & J & ~\color{red}{J} & ~\color{red}{Z} & ~\color{red}{J} & \color{red}{Z} & ~J & ~Z & ~Z & J & $A2$  \\
     \hline
  \begin{array}{@{}c@{}@{}c@{}}
    \bf{01} \\
    \bf{01} \\
  \end{array}
      & ~\color{red}{J} & \color{red}{J} & ~\color{red}{Z} & \color{red}{Z} & ~J & J & ~K & J & ~Z & ~J & ~J & Z & ~J & ~Z & ~J & Z & $C4$ \\
          \begin{array}{@{}c@{}@{}c@{}}
    \bf{11} \\
    \bf{01} \\
  \end{array}
       & ~\color{red}{J} & \color{red}{J} & ~\color{red}{Z} & \color{red}{Z} & ~J & J & ~J & K & ~J & ~Z & ~Z & J & ~Z & ~J & ~Z & J & $D3$    \\
       \hline
            \begin{array}{@{}c@{}@{}c@{}}
    \bf{00} \\
    \bf{10} \\
  \end{array}
       & ~J & Z & ~\color{red}{Z} & \color{red}{J} & ~\color{red}{Z} & \color{red}{J} & ~Z & J & ~K & ~J & ~J & J & ~J & ~J & ~Z & Z & $E5$     \\
          \begin{array}{@{}c@{}@{}c@{}}
    \bf{10} \\
    \bf{10} \\
  \end{array}
       & ~Z & J & ~\color{red}{J} & \color{red}{Z} & ~\color{red}{J} & \color{red}{Z} & ~J & Z & ~J & ~K & ~J & J & ~J & ~J & ~Z & Z & $F6$   \\
   \begin{array}{@{}c@{}@{}c@{}}
    \bf{00} \\
    \bf{11} \\
  \end{array}
      & ~J & Z & ~\color{red}{J} & \color{red}{Z} & ~\color{red}{Z} & \color{red}{J} & ~J & Z & ~J & ~J & ~K & J & ~Z & ~Z & ~J & J & $G7$   \\
         \begin{array}{@{}c@{}@{}c@{}}
    \bf{10} \\
    \bf{11} \\
  \end{array}
     & ~Z & J & ~\color{red}{Z} & \color{red}{J} & ~\color{red}{J} & \color{red}{Z} & ~Z & J & ~J & ~J & ~J & K & ~Z & ~Z & ~J & J& $H8$     \\
     \hline
     \begin{array}{@{}c@{}@{}c@{}}
    \bf{01} \\
    \bf{10} \\
  \end{array}
    & ~J & Z & ~J & Z & ~Z & J & ~J & Z & ~J & ~J & ~Z & Z & ~K & ~J & ~J & J & $G8$     \\
            \begin{array}{@{}c@{}@{}c@{}}
    \bf{11} \\
    \bf{10} \\
  \end{array}
     & ~Z & J & ~Z & J & ~J & Z & ~Z & J & ~J & ~J & ~Z & Z & ~J & ~K & ~J & J & $H7$   \\
             \begin{array}{@{}c@{}@{}c@{}}
    \bf{01} \\
    \bf{11} \\
  \end{array}
     & ~Z & J & ~J & Z & ~J & Z & ~J & Z & ~Z & ~Z & ~J & J & ~J & ~J & ~K & J & $F5$   \\
              \begin{array}{@{}c@{}@{}c@{}}
    \bf{11} \\
    \bf{11} \\
  \end{array}
  & ~J & Z & ~Z & J & ~Z & J & ~Z & J & ~Z & ~Z & ~J & J & ~J & ~J & ~J & K & $E6$  \\
  \end{array}
$$
\caption*{Figure 8: The adjacency matrix, $B_{e,2}$, of the subgraph $\Delta$ of $\Gamma _{e,2}$}
\end{table}
Let $(n,k,\lambda,\mu)$ be the parameters of the affine polar graph $VO^+(2e,2)$ as a strongly regular graph.
We have to check that the obtained graph $\Gamma_{e,2}$ is a strictly Neumaier graph.
Note that the vertices
$$
\left(
  \begin{array}{ccc|cc}
   \ast &\ldots & \ast & * & 1 \\
    0 &\ldots & 0 & 1  & * \\
  \end{array}
\right)
$$
induce a $2^{e-1}$-regular $2^e$-clique in $\Gamma_{e,2}$ as well as in $\Gamma_{e}$.
Let us check that any pair of vertices in $\Gamma_{e,2}$ is OK,
i.e. any two adjacent vertices have $\lambda$ common neighbours.
Also, we investigate which values of $\mu$ occur in $\Gamma_{e,2}$.

Let us consider any two vertices inside of $\Delta$.
The notation in the right column of the matrix in Figure 8 means the following.
Two block-rows have the same letter if and only if any row from one block-row
 and any row from the other block-row correspond to non-adjacent
vertices having $\mu+2^{e-1}$ common neighbours; two block-rows have the same number if and only if any row from one block-row
 and any row from the other block-row correspond to non-adjacent
vertices having $\mu-2^{e-1}$ common neighbours. Otherwise, every two non-adjacent vertices corresponding to rows of this submatrix have $\mu$ common neighbours. Any two adjacent vertices have $\lambda$ common neighbours. This means that
all pairs of vertices inside of $\Delta$ are OK.

Let us consider any two vertices outside of $\Delta$.
Their neighbours and, consequently, their common neighbours are preserved by the switching.
This means that all pairs of vertices outside of $\Delta$ are OK.

Let us consider a vertex $x$ in $\Delta$ and a vertex $y$ outside of $\Delta$.
If the neighbours of $x$ are preserved by the switching, then $x$,$y$ are OK.
Assume that the neighbours of $x$ are switched. Then the vertices $x$,$y$ are OK
since the vertex $y$ is adjacent to half of the vertices of each block of $\Delta$.
In fact, the vertex $y$ is presented by a matrix
$$
\left(
  \begin{array}{ccc|cc}
    y_1 & \ldots &y_{2e-5} & y_{2e-3} & y_{2e-1} \\
    y_2  & \ldots & y_{2e-4} & y_{2e-2} & y_{2e} \\
  \end{array}
\right),
$$
where there is at least one non-zero among $y_2, y_4, \ldots, y_{2e-4}$.
Without losing of generality, assume that $y_2 = 1$.
Let us show that $y$ is adjacent to half the of vertices in a block
$$
\left(
  \begin{array}{ccc|cc}
   \ast &\ldots & \ast & a & b \\
    0 &\ldots & 0 & c  & d \\
  \end{array}
\right).
$$
We have
$$
y +
\left(
  \begin{array}{ccc|cc}
   \ast &\ldots & \ast & a & b \\
    0 &\ldots & 0 & c  & d \\
  \end{array}
\right)
=
\left(
  \begin{array}{ccc|cc}
   \ast &\ldots & \ast & a' & b' \\
    1 &\ldots & y_{2e-4} & c'  & d' \\
  \end{array}
\right)
=
$$
$$
=
\left(
  \begin{array}{ccc|cc}
   0 &\ldots & \ast & a' & b' \\
    1 &\ldots & y_{2e-4} & c'  & d' \\
  \end{array}
\right)
\bigcup
\left(
  \begin{array}{ccc|cc}
   1 &\ldots & \ast & a' & b' \\
    1 &\ldots & y_{2e-4} & c'  & d' \\
  \end{array}
\right) = Y_0 \cup Y_1.
$$
Note that $|Y_0| = |Y_1|$, and the form $Q$ has value $0$ on one of the sets $Y_0,Y_1$ and value $1$ on the other.
We have proved that the switching preserves the number of common neighbours $x$ and $y$, completing the proof of the theorem. $\square$

\medskip
\section{Concluding remarks}
\medskip

There are four known non-isomorphic strictly Neumaier graphs with parameters $(24,8,2;1,4)$, all of which are vertex-transitive. An interesting open problem is to determine all strictly Neumaier graphs with these parameters (up to isomorphism). This will complete the classification of strictly Neumaier graphs on at most $24$ vertices.	

For $e=3$ and $4$, the two generalisations in Section 5 are known to give non-isomorphic graphs. We conjecture that the $i^{th}$ element of the first sequence of graphs is not isomorphic to the $i^{th}$ element of the second sequence of graphs, except for the value $i=1$.

Both of the constructions in Section 5 involve taking two pairs of disjoint regular cliques, and carrying out a switching between the cliques in each pair. Starting with the graph $VO^{+}(6,2)$, it can be shown computationally that any two such consecutive switchings between regular cliques give rise to only two distinct strictly Neumaier graphs, each of which appear in one of the above constructions. We also note that we can continue to apply switchings on disjoint regular cliques, and obtain many new strictly Neumaier graphs with the same parameters. For example, in this way we can show that there are at least $4$ non-isomorphic strictly Neumaier graphs with the same parameters as $VO^{+}(6,2)$. A natural question to ask is how many non-isomorphic strictly Neumaier graphs can we construct in this manner. We hope to use this iterative process to observe prolific constructions of strictly Neumaier graphs, similar to some prolific constructions of strongly regular graphs (see Wallis \cite{W71}, Fon-Der-Flaass \cite{F02}, Cameron $\&$ Stark \cite{CS02} and Muzychuk \cite{M07}).

The above constructions show that the nexus of a clique in a strictly Neumaier graph is not bounded above by some constant number. However, all known Neumaier graphs contain regular cliques with nexus $2^{j}$, for $j$ a non-negative integer. So we ask if there exist strictly Neumaier graphs containing regular cliques with nexus not a power of two? Finally we ask if we can generalise the above constructions to the case $q$ an arbitrary prime power, which would give strictly Neumaier graphs containing a regular clique with nexus a prime power.

\section*{Acknowledgments} \label{sec3}

We would like to express our gratitude to Leonard Soicher and Alexander Gavrilyuk for introducing the authors and their continued support. We are also grateful to Jack Koolen and Gary Greaves for their advice and suggestions on the topics discussed. Finally, we would like to thank Derek Holt and Gordon Royle for providing us with their enumeration of small vertex-transitive edge-regular graphs.

%% Not needed for arXiv
%%\section*{References}

\end{document}